\documentclass[]{article}

\usepackage[T1]{fontenc}
\usepackage{csquotes}
\usepackage[english]{babel}
\usepackage[a4paper]{geometry}

\usepackage{orcidlink}
 
\usepackage{graphicx}

\usepackage{amsfonts, amsmath, amsthm, amssymb}

\usepackage{enumitem}

\usepackage{subcaption}  

\usepackage{algorithm}
\usepackage[italicComments=true]{algpseudocodex} 

\usepackage{standalone}

\usepackage{xcolor}
\usepackage{pgf}
\usepackage{tikz}
\usepackage{pgfplots}  
\pgfplotsset{compat=1.15}
\usetikzlibrary{plotmarks,backgrounds,positioning}
\usetikzlibrary{calc}
\usetikzlibrary{patterns}
\usetikzlibrary{patterns.meta}

 \usepackage{hyperref}

\usepackage
[
maxnames=5, 
giveninits=true, 
isbn=false,
backref=true
]
{biblatex}

\DeclareSourcemap{
  \maps[datatype=bibtex]{
    \map[overwrite]{
      \step[fieldsource=doi, final]
      \step[fieldset=url, null]
      \step[fieldset=eprint, null]
    }  
  }
}

\usepackage{comment}

\usepackage{todonotes}

\definecolor{amber}{rgb}{1.0, 0.75, 0.0}
\definecolor{apricot}{rgb}{0.98, 0.81, 0.69}
\definecolor{aquamarine}{rgb}{0.5, 1.0, 0.83}

\newcommand{\mh}[1]{{#1}}
\newcommand{\ck}[1]{{#1}}
\newcommand{\ph}[1]{{#1}}

\theoremstyle{plain}
\newtheorem{theorem}{Theorem}

\theoremstyle{remark}
\newtheorem{remark}[theorem]{Remark}

\theoremstyle{definition}

\DeclareMathOperator{\divergence}{div}
\DeclareMathOperator{\trace}{tr}
\newcommand{\dx}{{\ \mathrm{d}x}}
\newcommand{\ds}{{\ \mathrm{d}s}}
\newcommand{\R}{\mathbb{R}}
\newcommand{\aphi}{{a}}
\newcommand{\bphi}{{b}}
\newcommand{\cphi}{{c}}
\newcommand{\id}{{\rm id}}

\DeclareMathOperator{\diam}{diam}

\begin{document}

\title{Combining diffuse and sharp interface methods in shape optimisation}

\author{Philip J.~Herbert%
\footnote{Department of Mathematics,
    University of Sussex, United Kingdom} %
\orcidlink{0000-0002-6513-1728},\,%
Michael Hinze%
\footnote{Mathematical Institute, University of Koblenz, Germany} %
\orcidlink{0000-0001-9688-0150},\,%
and %
Christian Kahle%
\footnotemark[2]
\orcidlink{0000-0002-3514-5512}\,%
}

\maketitle

\begin{abstract}
    We develop a concept for the numerical treatment of shape optimization problems based on the combination of phase field and sharp interface methods. On the one hand, phase field methods are very well suited to numerically determine the shape, size and topology of a sought domain, but on the other hand they have problems to sharpen out domains where they e.g. should develop corners. However, this is the strength of a sharp-interface approach developed in our group, which provides shape updates in the Lipschitz topology. This leads to a two-stage process that first determines an optimized shape using the phase field method. The resulting domain is the starting solution for the sharp interface shape optimization method. 
    Both methods are  discretized with the finite element method.  
    The starting mesh for the sharp method is constructed from the finite element mesh of the optimal phase field solution using its properly post processed  zero-level set. We describe this construction process in detail and investigate the performance of our method on a selection of test problems from the literature and from applications.
\end{abstract}

\section{Contribution}
We propose a two-step numerical procedure for PDE constrained shape optimization. In the first step we use numerical phase field shape and topology optimization to
\begin{enumerate}[label=\alph*)]
\item determine the topology of the sought object, 
\item to adaptively construct an aligned mesh for the sought shape(s), and to
\item provide triangulated initial shapes for the sharp interface approach. \label{en:contr:iii}
\end{enumerate}
The most technically demanding part is step~\ref{en:contr:iii}. 
We use the zero level set of the optimal phase field as shape boundary of the initial domain for the sharp interface approach. This defines a polygonal domain. We then construct a triangulation of this polygonal domain from the triangulation of the optimal phase field. In the second step, we use this triangulation as the initial configuration for the finite element method in the sharp interface approach.

\begin{figure}
    \centering
    \includegraphics[trim= 0 200 0 0, clip, width=0.45\linewidth]{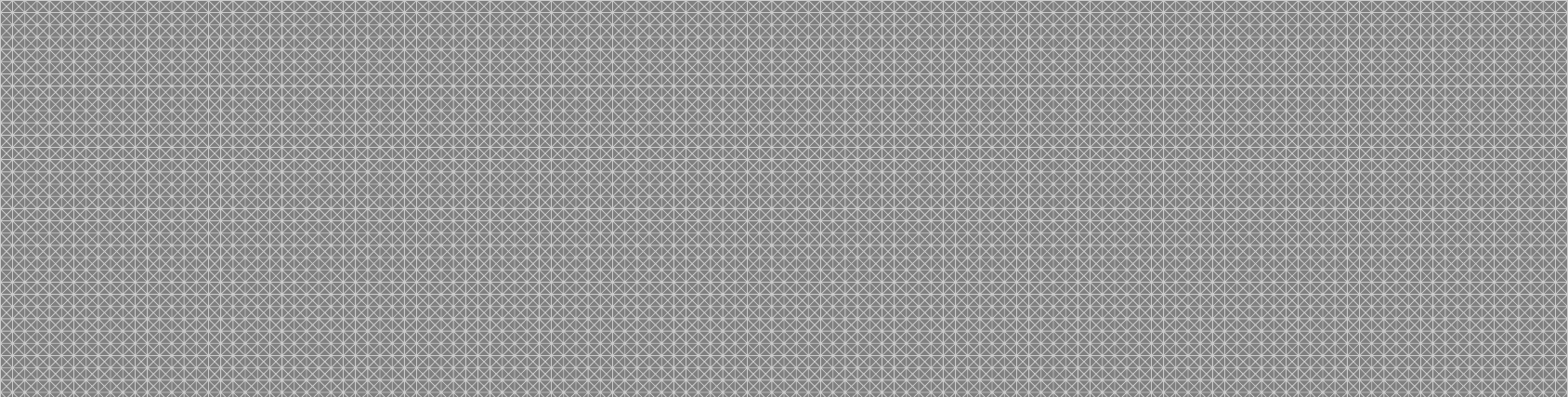}%
    \hspace{1ex}%
    \includegraphics[trim= 0 200 0 0, clip, width=0.45\linewidth]{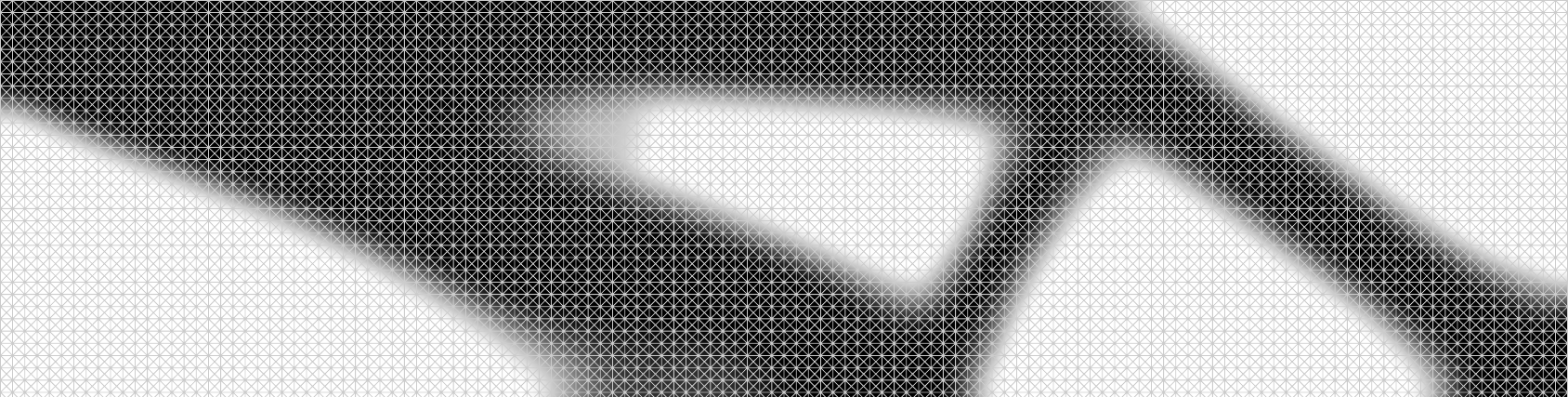}\\[1ex]
    \includegraphics[trim= 0 0 0 200, clip,width=0.45\linewidth]{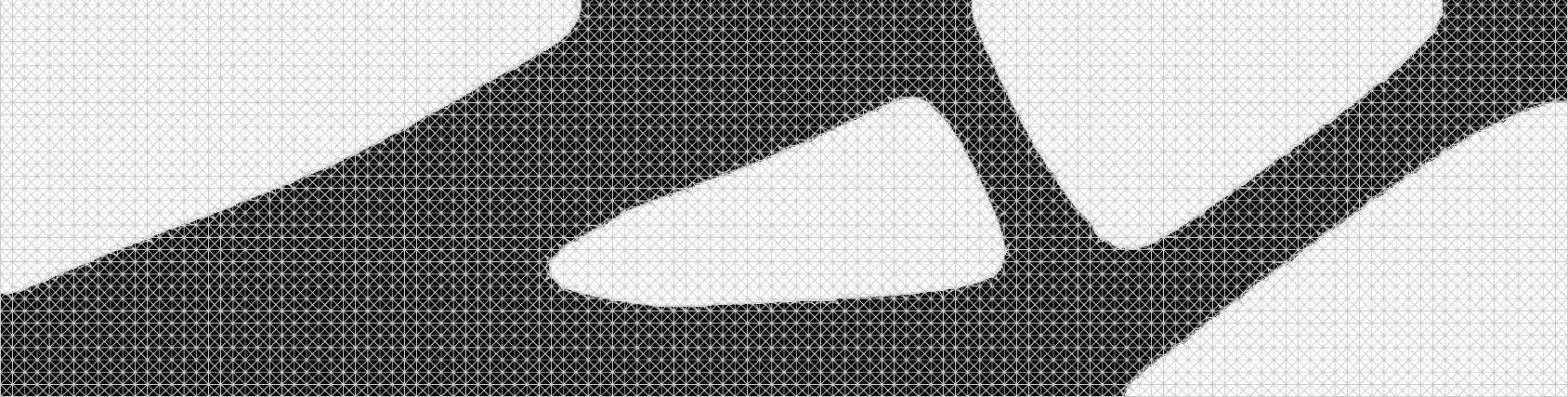}%
    \hspace{1ex}%
    \includegraphics[trim= 0 0 0 200, clip,width=0.45\linewidth]{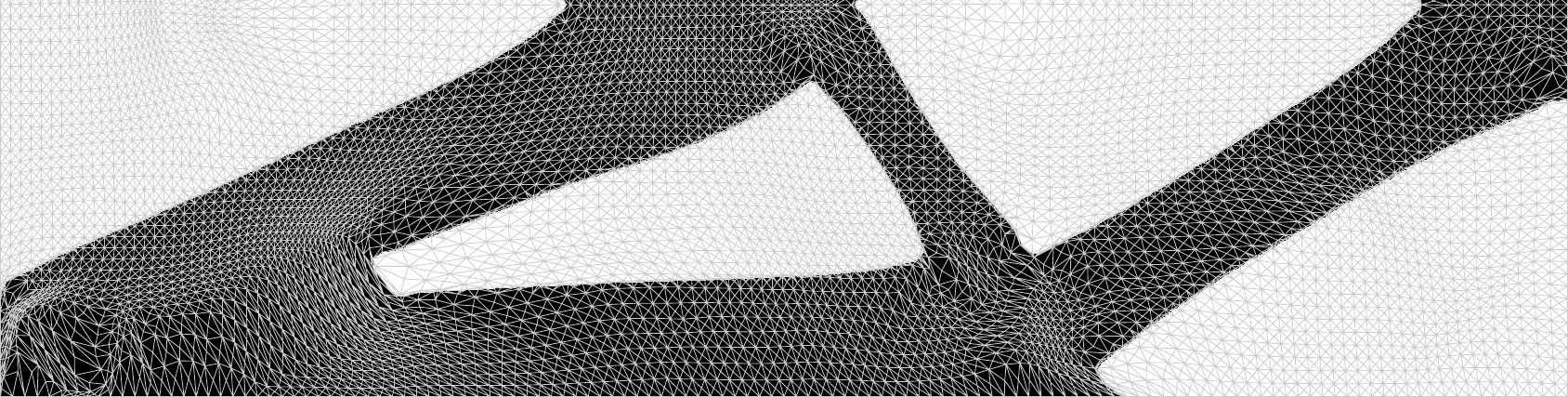}
        
    \caption{Production line: Phase-field initialization (top left), phase-field optimized (top right), sharp initialized (bottom left), sharp optimized (bottom right) for structural topology optimization. Since the result is symmetric with respect to the horizontal axis we only show one half of the domain.
    For the details we refer to Section~\ref{sec:Elasticity}.
    }
    \label{fig:pl}
\end{figure}

Our approach involves reusing as much numerical infrastructure as possible from the phase field approach in the sharp interface setting. This procedure is summarised in the production line, which is demonstrated in Figure~\ref{fig:pl} by the example detailed in Section \ref{sec:Elasticity}.

\section{Introduction}

In this work, we focus on the practical implementation of PDE-constrained shape optimisation problems where the topology of a minimiser is unknown or where only limited information about the size and shape of a possible minimiser is available. The PDE-constrained shape optimisation problem takes the following form
\begin{equation}
\tag{\ensuremath{P}}
\label{eq:sharp-P}
    \min_{\Omega ,y} J(\Omega, y) 
    \quad \mathrm{s.t} \quad 
    e(\Omega,y) = 0 \text{ and }\Omega \in  \mathcal{S}_{ad}.
\end{equation}
Here, $\mathcal{S}_{ad}$ denotes a collection of admissible domains contained within the convex and bounded hold-all domain $D$.
 The PDE constraint is encoded by $e(\Omega,\cdot) = 0$, and $J$ denotes the cost functional (sometimes called the energy) which depends on the domain $\Omega$ and the solution $y$ of the PDE  $e(\Omega,y) = 0$. 
A collection of application settings is provided in the numerical examples Section~\ref{numerical-examples}.

The task of PDE constrained shape optimisation is to find a domain that minimises an energy which depends on the solution to a PDE inside the domain.
Such problems arise in many practical applications, including fluid dynamics and elasticity, for example. 
In practice, mesh-based methods are used to solve the PDE. 
As the shape evolves during the optimisation process, the mesh is often moved simultaneously with the domain to avoid the need for remeshing. For the numerical solution of the PDE, it is crucial that the quality of the mesh is maintained as much as possible, and that the mesh elements do not overlap. 
Many methods have been developed to ensure that mesh quality does not deteriorate significantly, see e.g.~\cite{blauth2026enforcing,chen2026inertial,EltHerLoa20,iglesias2018two,luft2021simultaneous,LufSch21-A} for Hilbertian shape optimization approaches. 
Alternatively, Banach space methods have recently proven to be very promising in this context, also allowing convergence analyses of the underlying numerical algorithms \cite{DecHerHin22,DecHerHin24-App, HerPinSie23, DecHerHin25-NA,MulKuhSie21,MulPinRun23,starke2024shape}.
For further references of (mostly Hilbertian) shape optimization methods we refer to the excellent overview article \cite{AllDapJou21}.

However, it is noteworthy that methods which modify the mesh are restricted to shapes with a fixed topology.
In particular, if the initial shape is incorrect, this can significantly degrade the quality of the mesh as the optimisation procedure attempts to correct it.

\begin{figure}
    \centering
    \includegraphics[width=0.45\linewidth]{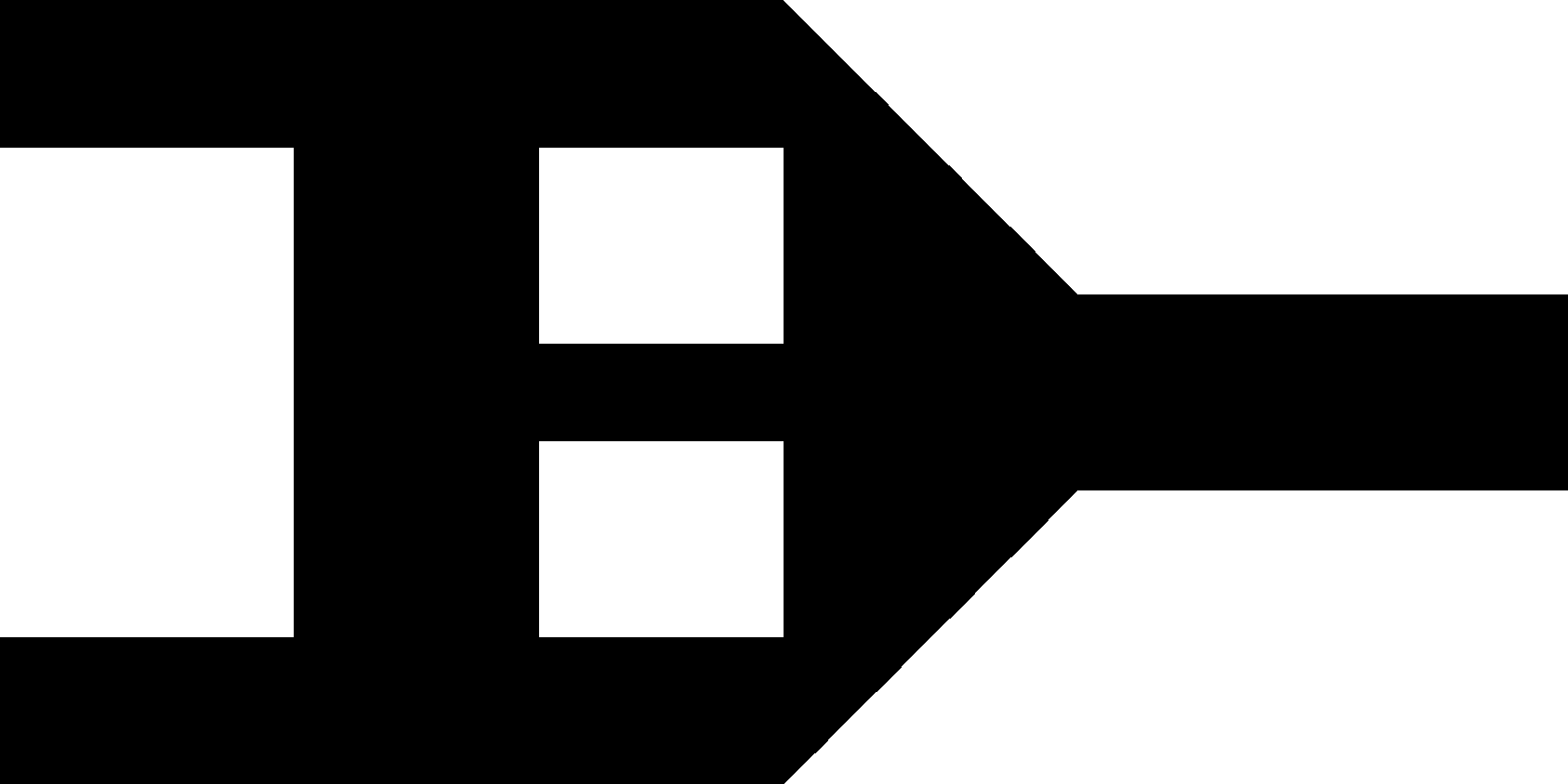}
    \hfill
    \includegraphics[width=0.45\linewidth]{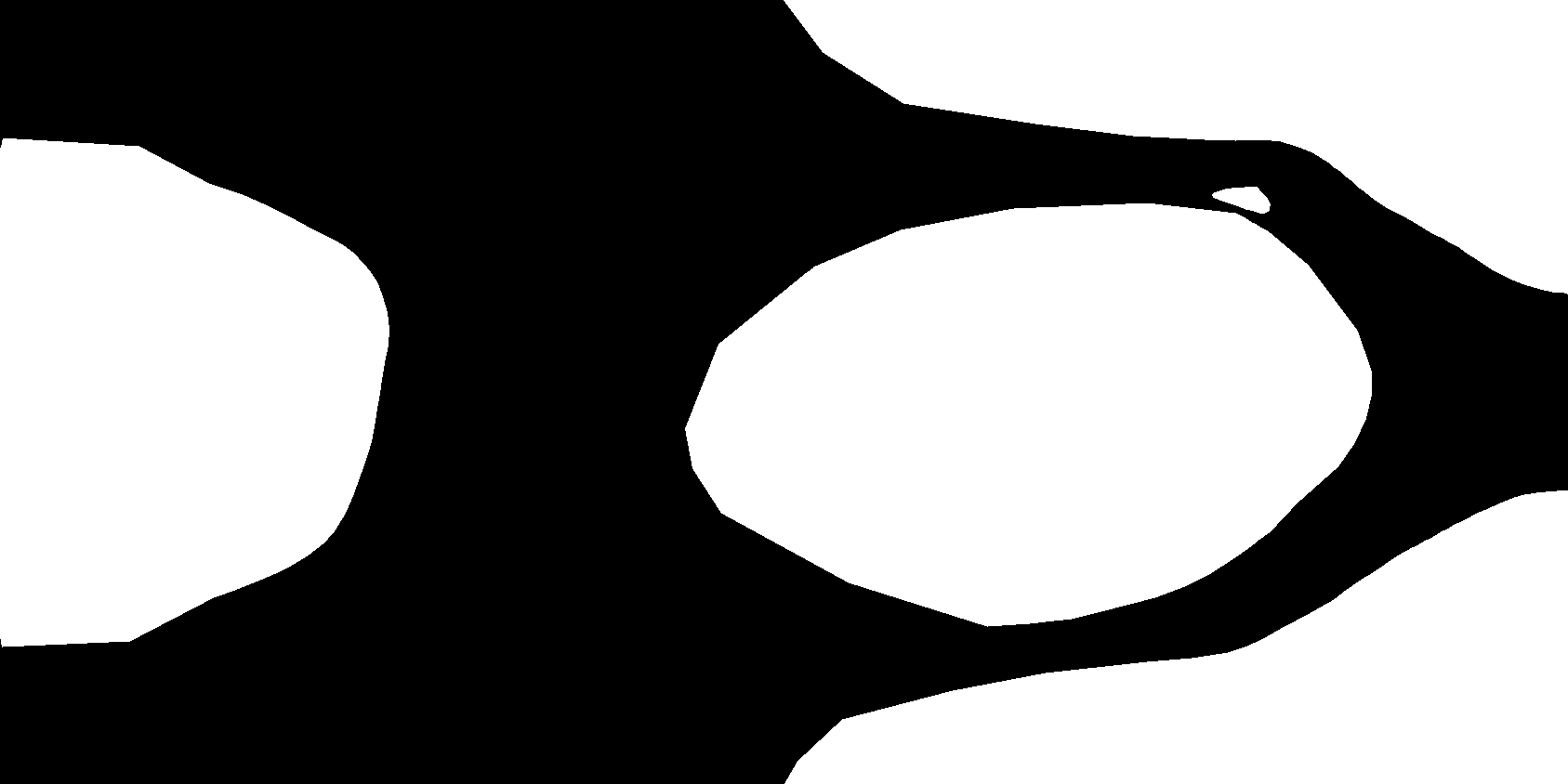}
    \caption{Shape optimisation of a problem in linear elasticity.
    The image on the left is an initial guess, the image on the right is the optimisation after 40 steps of a $W^{1,\infty}$ steepest descent method with volume preservation.
    We can see from this experiment, that we have gone from two equal sized holes to one large hole and one much smaller hole.
    One may interpret that as the optimum having a single, larger hole.
    Such a large change in the shape has caused a significant change in mesh quality, with both holes having the same number of edges along their boundaries.}
    \label{fig:exampleTopOpt}
\end{figure}

\ck{
This is illustrated by the images in Figure~\ref{fig:exampleTopOpt}. Here, we consider an example of structural optimisation. A beam is fixed on the left of the domain and a force is applied on the right. On the left, we show the initial domain for the shape optimisation process, and on the right, we show the result after 40 steps of a  $W^{1,\infty}$ steepest descent method. 
It can be seen that one of the two holes on the right has been enlarged while the other has shrunk.
It is clear, therefore, that the initial guess does not provide the correct topology of the optimised domain, and that the optimisation approach cannot change the topology of the domain. Since the shape is represented using a mesh, both holes in the initial guess are represented by the same number of triangles. This indicates significant mesh degradation, as both the tiny and the large hole in the optimised shape are represented by the same number of triangles.
}

\medskip

To alleviate the challenges of an incorrect topology or large changes in the domain, we propose to use a phase field method to generate an initial guess for the shape.
Phase fields have been a frequent tool used in shape optimisation in recent years, see e.g. \cite{
BlankEtAl-2012-PhaseFieldStructuralTopoOpt,%
GarHHk-2015-NumApproxPhaseTopoOptFluids,%
GarckeHKKL-2023-PhaseFieldShapeOptLaplace,%
GarckeLamNuernbergSignori-2023-OverhangPenalizationManufacturingTopoOpt,%
GarckeHKK-2024-MultiPhaseElasticSpectral,%
}
\cite{%
Evgrafov-2005-LimitsOfPorousMediumTopologyOptimization%
,ZhouWang-2006-MulitmaterialStructuralTopoOptGeneralizedCHMultiphase%
,Bourdin-Chambolle%
,Burger-Stainko%
,TakezawaNishiwakiKitamura-2010-ShapeTopoOptPhaseFieldSensitivity%
,Penzler%
,Wallin2012%
,DedeBordenHughes-2012-IsogeometricalTopologyOpt%
,Bartels-2015-RobustnessErrorEstimatesPhaseFieldTopoChange%
,Carraturo%
,Auricchio2%
,AlmiStefanelli%
}.

This approach transforms the shape optimisation problem into a PDE-constrained optimisation problem, where the control is given by a phase-field parameter acting as coefficient in the underlying PDEs. This enables the entire PDE-constrained optimisation framework to be invoked for analytical and numerical treatment  \cite{HPUU2009}.
This approach has many advantages, such as the availability of monolithic solvers that can be implemented straightforwardly, the fact that an initial domain need not be specified, and the topology need not be fixed. However, the phase field method has many parameters that need to be adjusted for proper implementation and has limited ability to resolve critical geometric quantities.

In the present work we combine the advantages of the phase-field and sharp-interface approaches to PDEs in constrained shape optimisation. We demonstrate the practicality of this combined approach using selected examples from the literature.

Finally, we would like to mention that other methods combining topology and shape optimization are also proposed in the literature. 
The paper \cite{wang2022sequentially} proposes a sequential coupling of gradient-based topology and shape optimization. The authors in \cite{auricchio2024analysis} analyze a combination of filter and phase field methods. In \cite{wang2007extended}, an approach is presented that couples the level set method with sharp shape optimization.

\mh{Let us also note that we are aware of other approaches which might be used to construct a starting situation for the sharp interface method, as e.g.~the unfitted finite element method \cite{barrett1987fitted}, where then the optimal phase field grid is used together with the zero level set of the phase field function to initialize the sharp interface method.}

\paragraph{Outline.} 
In Section~\ref{sec:problemsGeneral} we introduce the general shape and topology optimization problem that we investigate in the sharp setting and introduce how corresponding phase field approximations can be formulated in Section~\ref{sec:problemsGeneralPF}.
The coupling from phase field solution to sharp interface problem is described in Section~\ref{sec:Tech} and finally in Section~\ref{sec:Poisson}, Section~\ref{sec:Elasticity}, and Section~\ref{sec:NavierStokes} we apply the concept numerically to three prototypically physical systems, namely the Poisson equation,  linear elasticity, and  the Navier--Stokes equation.


\section{PDE constrained shape optimisation and its approximation by phase fields}
\label{sec:problemsGeneral}

In this section, we present the shape optimization methods used in this work to address task \eqref{eq:sharp-P} in more detail  \ck{and show the general concept of phase field approximations to these shape optimization problems.}

For convenience, we repeat problem~\eqref{eq:sharp-P} here
\begin{equation}
\tag{\ensuremath{P}}
    \min_{\Omega ,y} J(\Omega, y) 
    \quad \mathrm{s.t} \quad 
    e(\Omega,y) = 0 \text{ and }\Omega \in  \mathcal{S}_{ad}.
\end{equation}
Let us assume that the PDE $e(\Omega,y)=0$ admits a unique solution $y(\Omega)$. Then we can define the \textit{shape functional} $\mathcal J(\Omega):= J(\Omega,y(\Omega))$ and the shape optimization problem \eqref{eq:sharp-P} can be equivalently rewritten in the reduced form 
\begin{align*}
    \min_{\Omega \in \mathcal S_{ad}} \mathcal J(\Omega).
\end{align*}

\subsection{Method of mappings}
\label{sec:probGen-sharp}
To constructively find a (local) minimiser - or stationary point - of a shape functional, we use the steepest descent method. Given a starting domain/shape $\Omega_0$, in the $n$-th step of this method the shape $\Omega_n$ is updated according to $\Omega_{n+1} = (\id + t_n V_n)(\Omega_n)$, where $t_n\in (0,1)$ is some step size and $V_n\colon \Omega_n \to \mathbb{R}^d$ is a suitable descent vector field which guarantees $\mathcal J(\Omega_{n+1}) < \mathcal J(\Omega_n)$, i.e.~that the shape functional decreases. Here, we use the method proposed in \cite{DecHerHin24-App,DecHerHin25-NA} and determine the vector field $V_n$ as the direction of steepest descent for the shape functional in the $W^{1,\infty}$ topology, where $t_n \in (0,1)$ in the steepest descent method is determined by means of the Armijo--Goldstein rule. This guarantees that the mapping $\id + t_n V_n$ is a bi-Lipschitz map and among other things ensures that the topology of the domains $(\Omega_n)_n$ remains unchanged throughout the iterations, and also that open domains are mapped to open domains. In particular, the latter  ensures that  e.g.~the weak formulation of the PDEs with homogeneous Dirichlet boundary data is well defined on every domain $\Omega_n$.

Some of the problems to be considered here involve geometric constraints.
One simple example of such a constraint is that the volume is fixed.
This alters the way in which the vector field $V_n$ is selected.
This issue is addressed in \cite{MulPinRun23} for a $p$-Laplace regularised version of the $W^{1,\infty}$ method of \cite{DecHerHin25-NA} and for the full $W^{1,\infty}$ problem, was considered in \cite{HerPinSie23}.
While convergence results as in \cite{DecHerHin25-NA} are not available for this geometrically constrained problem, they are of interest and will be the subject of future work.

\bigskip

\mh{To discretize we use the finite element method as proposed in \cite{DecHerHin24-App,DecHerHin25-NA}. For simplicity, assume that the initial shape $\Omega_0$ is a polygon and admits a regular triangle mesh $\mathcal{T}_0$. With the domain and the mesh we associate the finite element space consisting of piecewise linear and continuous functions. The meshes together with the finite element spaces are then updated as the domains are with
\begin{equation}
    \mathcal{T}_{n+1}
    =
    \{ (\id + t_n V_n)(T) : T \in \mathcal{T}_n\},
\end{equation}
where now $V_n$ denote piecewise linear and continuous vectorfields defined on the triangulation $\mathcal T_n$, and which are chosen as steepest descent directions in the $W^{1,\infty}$ topology of the discretized shape  functional $\mathcal J_h(\Omega) := J(\Omega,y_h(\Omega))$ with $y_h(\Omega)$ denoting the finite element approximation of the PDE solution $y(\Omega)$. Since the respective mappings $\id +t_nV_n$ are bi-Lipschitz for $t_n \in (0,1)$, all triangulations $\mathcal T_n$ are regular.}

Since the updates which we consider are bi-Lipschitz, hence homeomorphic, topological quantities are preserved.
However, this  is problematic 
if the topology of the sought domain is not known in advance and respected by $\Omega_0$.

Even if the initial guess of $\Omega_0$ has the correct topology for the optimal shape, significant deformation of the original  shape may be necessary to achieve it.
Since the mesh moves with the shape, this can result in a significant reduction in mesh quality and potentially lead to an inaccurate approximation of the underlying PDE.

These two issues can be addressed by having a good guess for the initial shape.
Here, we propose utilising a phase field approach to achieve this.

\begin{remark}
    \mh{In general, the existence of solutions to shape optimization problems with PDE constraints is by no way clear. In e.g. \cite[Chapter 4.2]{henrot2018shape} a couple of PDE constrained shape optimization problems for the Poisson problem are discussed, which do not admit a solution. Existence results depend strongly on the topology imposed on admissible domains, geometric compactness assumptions, and stability properties of the PDE under domain perturbations.
    We note that \cite[Chapter 5]{BucBut05} discusses some of the settings wherein optimal domains can be shown to exist. It is an advantage of the phase field approximation of PDE constrained shape optimization problems, that among other things existence of solutions can be shown, and also that the Hilbert space framework gives more flexibility in deriving analytical results than the sharp formulation of the problem, see Remark~\ref{rem:PF:exsol}}.
    
\end{remark}

\subsection{Approximation with a phase field approach}
\label{sec:problemsGeneralPF}

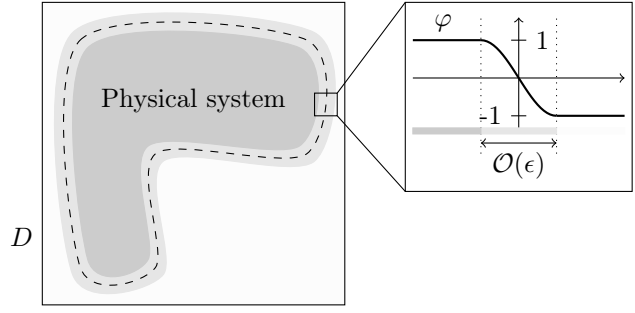
\begin{figure}
    \centering
     
  \begin{tikzpicture}[scale=2]

  \colorlet{clOmP}{black!20}
  \colorlet{clOmM}{black!01}

\draw[fill,color=clOmM] (0,0) rectangle(2,2);	
\draw[color=black] (0,0) rectangle(2,2);	
\draw (0.0,0.45) node[anchor=east] {$D$};

\draw[fill,color=black!10] plot[smooth cycle] coordinates{(0.2,0.2)(0.20,1.8)(1.8,1.8)(1.8,0.95)(0.8,0.95)(0.8,0.2)};

\draw[dashed] plot[smooth cycle] coordinates{(0.25,0.25)(0.25,1.75)(1.75,1.75)(1.75,1)(0.75,1.00)(0.75,0.25)};

\draw[fill,color=clOmP] plot[smooth cycle] coordinates{(0.3,0.3)(0.3,1.7)(1.7,1.7)(1.7,1.05)(0.7,1.05)(0.7,0.3)};


\draw (1.0,1.35) node[anchor=center] {Physical system};

\draw (1.8,1.25) rectangle ++(0.15,0.15);

\draw (1.95,1.25) -- (2.4,0.75);
\draw (1.95,1.4) -- (2.4,2.0);

\begin{scope}[xshift=2.4cm,yshift=0.75cm]

\draw (0,0) rectangle(1.5,1.25);	  

\draw[->] (0.05,0.75) -- (1.45,0.75); 
\draw[->] (0.75,0.40) -- (0.75,1.15); 

\draw (0.70,1.0)--(0.80,1.0);
\draw (0.70,0.5)--(0.80,0.5);
\draw (0.8,1.0) node[anchor=west] {1};
\draw (0.7,0.5) node[anchor=east] {-1};

\draw[dotted] (0.5,0.25) -- (0.5,1.15); 
\draw[dotted] (1.0,0.25) -- (1.0,1.15);

\colorlet{clphi}{black!100}
\draw[thick,color=clphi,domain=-1:1] plot (0.75 + \x/4,{0.75-0.25*(sin(pi/2*\x r))});
\draw[thick,color=clphi] (0.05,1.0) -- (0.5,1.0);
\draw[thick,color=clphi] (1.0,0.5) -- (1.45,0.5);
\draw (0.25,1.0) node[anchor=south]{\color{clphi}$\varphi$};

\draw[fill,color=clOmM] (1.0,0.38) rectangle ++(0.45,0.04);
\draw[fill,color=black!10] (0.5,0.38) rectangle ++(0.5,0.04);
\draw[fill,color=clOmP](0.05,0.38) rectangle ++(0.45,0.04);
\draw[<->] (0.5,0.32)--(1.0,0.32);
\draw (0.75,0.32) node[anchor=north] {$\mathcal O (\epsilon)$};

\end{scope}

\end{tikzpicture}
    
    \caption{Sketch of the phase field approximation. 
     The values of $\varphi\in[-1,1]$ are indicated in gray
      and the zero level set of $\varphi$ is indicated by the dashed curve. 
     The sign of $\varphi$ splits $D$ into two subdomains.   
    We identify $\Omega = \{ \varphi \geq 0\}$. 
    Inside $\Omega$ a physical system governed by a PDE is satisfied.
     }
    \label{fig:PF-sketch}
\end{figure}

In the following we summarise how the shape optimisation problem \eqref{eq:sharp-P} can be properly approximated by a phase field approach.
In the phase field  approach, admissible domains $\Omega$ are represented by a \emph{smooth indicator function} $\varphi$, called phase field, satisfying $|\varphi|\leq 1$, 
according to $\Omega = \{\varphi \ge 0\}$, compare Figure~\ref{fig:PF-sketch}.
\ck{We construct the phase field $\varphi$ such that it takes values $\varphi\in\{\pm1\}$
in most parts of the hold all domain $D$, while indicating $\Omega = \{\varphi \geq 0\}$ and $D\setminus \Omega = \{\varphi<0\}$.
Its value changes smoothly between these subdomains over a length scale of order $\epsilon>0$ for some small value of $\epsilon$. Since $|\varphi| < 1$ in this region, this approach is also referred to as  \textit{diffuse}.
}
We refer e.g.~to \cite{BlankEtAl-2012-PhaseFieldStructuralTopoOpt,GarHHk-2015-NumApproxPhaseTopoOptFluids,GarckeHKKL-2023-PhaseFieldShapeOptLaplace} for a detailed description of the approach.

In our situation we replace the shape optimization problem \eqref{eq:sharp-P} by a phase field approximation \eqref{eq:PF-Peps}.
On a given hold all domain $D$, we  introduce a phase field $\varphi \in H^1(D) \cap L^\infty(D)$, $|\varphi|\leq 1$ and consider the minimization problem
\begin{equation}
\tag{\ensuremath{P^\epsilon}}
\label{eq:PF-Peps}
    \begin{aligned}
    \min_{\varphi,y^\epsilon}~& J^\epsilon(\varphi,y^\epsilon) 
    + \gamma E^\epsilon(\varphi)\\
    \mathrm{such~that}&\\
    &e^\epsilon(\varphi,y^\epsilon) = 0,\\
    &\varphi \in  \Phi_{ad}.
    \end{aligned}
\end{equation}
Here, $J^\epsilon$ and $e^\epsilon$ denote suitable  problem-dependent approximations of $J$ and $e$, respectively. For details we refer to the application section; Section~\ref{numerical-examples}.

Moreover, $ \gamma E^\epsilon(\varphi)$ denotes a regularization term, with a regularization parameter $\gamma>0$.
Here we choose the Ginzburg--Landau energy
\begin{align}
    E^\epsilon (\varphi) := \frac{1}{c_0}\int_D \left( \frac{\epsilon}{2} |\nabla \varphi|^2 + \frac{1}{\epsilon}\Psi(\varphi) \right) \dx
    \quad \mathrm{where} \quad 
    \Psi(\varphi) = 
    \begin{cases} 
    1-\varphi^2 & \mathrm{if }\ |\varphi|\leq 1, \\ 
    \infty & \mathrm{else},
    \end{cases}
\end{align}
which represents an approximation of the surface volume of the set $\{\varphi\geq 0\}$, see \cite{Modica}. The constant $c_0$ depends on the form of $\Psi$ and in our case is given by $c_0 = \frac{\pi}{2}$.
Due to its non-smooth nature, minimisation of $\Psi$ inherently contains the pointwise constraint $|\varphi| \leq 1$ on $\varphi$ and therefore gives rise to variational inequalities.
The geometric constraints are encoded in the set of admissible functions
$\Phi_{ad} \subset \{ \varphi \in H^1(\Omega) \mid |\varphi|\leq 1\} \subset H^1(D)\cap L^\infty(D)$. An explanation is provided in the applications section; Section~\ref{numerical-examples}.
Especially the constraint $|\varphi|\leq 1$ is treated in this way.

For the numerical realisation, we use an admissible triangulation to  discretise $D$
and piecewise linear, globally continuous finite elements to discretise $\varphi$.  
The state $y^\epsilon$ is then discretised according to the specific state equation as outlined in the examples below. 
 
\ck{
\begin{remark}
\label{rem:PF:exsol}
The phase field approximation formulates the shape and topology optimisation problem as a PDE constrained optimization problem with controls in the coefficients, where the control is given by the phase field function. This enables the use of well-established analytical concepts for this type of problem, see e.g.~\cite{HPUU2009}. 
In particular, we can typically expect to obtain existence results for solutions to the phase field approximation, see e.g.~\cite{GarHKL-2018-TopOptPF-StateConstraints,
GarckeHKK-2024-MultiPhaseElasticSpectral,
GarckeHKKL-2023-PhaseFieldShapeOptLaplace}.

We also stress that, in general, the topology of an optimised shape is not known a priori. Thus, topology changes may occur during the optimisation procedure. Methods that explicitly represent the shape may therefore break down in this situation, whereas the phase field approximation does not explicitly represent the topology of the shape, allowing topology changes to occur naturally.
\end{remark}
}

\begin{remark}
We would like to emphasise at this point that equation $\eqref{eq:PF-Peps}$  contains the Ginzburg--Landau energy $E^\epsilon(\varphi)$ as a regulariser to achieve minimisers in the form of a phase field.
The  term $E^\varphi(\varphi)$ converges  to the volume of $\partial \Omega$ for  $\epsilon \to 0$, see \cite{Modica}.
Therefore, any solution of $\eqref{eq:PF-Peps}$ inherently contains surface regularisation.
Whether or not this regularisation is also used for the sharp formulation \eqref{eq:sharp-P} depends on the problem and may affect the quality of the phase field solution as an initialiser for the sharp solution strategy.
\end{remark}



\section{Technical interface between sharp problem and phase field problem}
\label{sec:Tech}
The central aspect of our proposed method, which combines sharp interface and phase field methods for shape and topology optimisation, lies in generating a suitable initial domain, $\Omega_0$, for the sharp interface problem.
To retain mesh properties during sharp interface optimisation, $\Omega_0$ should be a good estimate of an optimal domain. In particular, the correct topology must be identified.
In this section we explain the generation of  $\Omega_0$ given a sufficiently optimal phase field.


Our strategy is detailed in two spatial dimensions in Algorithm~\ref{alg:p2s:phase2sharp}.
We currently do not consider three dimensions; however, such an approach could be devised, albeit it would be more technically complicated.

Consider the problem \eqref{eq:sharp-P} and its phase field approximation \eqref{eq:PF-Peps} on a hold-all domain $D$ with solution $\varphi$.
Let $\mathcal T_h$ denote a triangulation of $D$ and let $\varphi_h$ denote an approximation of $\varphi$ using piecewise linear, globally continuous finite elements.
Since $\varphi_h$ is piecewise linear and continuous, its zero level line is a polygon (where we assume that $\nabla \varphi_h \not= 0$ on $\{\varphi_h=0\}$) and the vertices of the polygon are located on edges of $\mathcal T_h$.
To find the initial shape, $\Omega_0$, for the sharp method, including a triangulation $\mathcal T_h^0$ of $D$ which represents $\Omega_0$ exactly, 
we define a first guess $\Omega_0^* = \{\varphi_h \geq 0\}$, i.e.~we consider the zero level-line of $\varphi_h$ as the boundary of $\Omega_0^*$.
This definition of $\Omega_0^*$ suggests, that based on $\mathcal T_h$ we in a first step generate a triangulation of $D$ that exactly represents $\Omega_0^*$ by adding the zero-level line of $\varphi_h$ to $\mathcal T_h$. This requires to add vertices on edges of $\mathcal T_h$ that intersect with the zero level line. However, this approach may result in triangles that are too thin if the intersection points are too close to existing vertices of the triangulation. $\mathcal T_h$. 
Therefore, newly added vertices occasionally need to be shifted towards the centre of their respective edges. This is realised in a second step where we post-process the triangulation to obtain a conforming triangulation of  $D$ which we denote by $\mathcal T_h^0$.  

\begin{figure}
    \centering
    {
    \tikzset{every picture/.style={scale=0.95}}
        \colorlet{clgr+}{black!20}
    \colorlet{clgr-}{black!10}

\begin{tikzpicture}[scale=4]
    
\begin{scope}

 	\coordinate (xP1) at (0.0,0.0);
    \coordinate (xP2) at (1.0,0.0);
    \coordinate (xP3) at (0.5,0.86);
    \coordinate (xP4) at (1.5,0.86); 
    \coordinate (xP5) at (2.0,0.0);


    \coordinate (xI13) at ($(xP1)!0.65!(xP3)$);
    \coordinate (xI23) at ($(xP2)!0.9!(xP3)$);
    \coordinate (xI24) at ($(xP2)!0.65!(xP4)$);
    \coordinate (xI45) at ($(xP4)!0.9!(xP5)$);

    \coordinate (xI13s) at (xI13);
    \coordinate (xI23s) at ($(xP2)!0.70!(xP3)$);
    \coordinate (xI24s) at (xI24);
    \coordinate (xI45s) at ($(xP4)!0.70!(xP5)$);

    \draw[fill,color=clgr-] (xP1)--(xP2)--(xI23)--(xI13)--cycle;
    \draw[fill,color=clgr-] (xP2)--(xI24)--(xI23)--cycle;
    \draw[fill,color=clgr-] (xI45)--(xI24)--(xP2)--(xP5)--cycle;


    \draw[fill,color=clgr+] (xI13) -- (xI23) -- (xP3) -- cycle;
    \draw[fill,color=clgr+] (xI23) -- (xP3) -- (xP4) -- (xI24) -- cycle;
    \draw[fill,color=clgr+] (xI24) -- (xI45) -- (xP4) -- cycle;
       

    \draw[very thick] (xP1) -- (xP2) -- (xP3) --cycle;
    \draw[very thick] (xP2) -- (xP4) -- (xP3) --cycle;
    \draw[very thick] (xP2) -- (xP4) -- (xP5) --cycle;

    \tikzset{
    dot/.style={circle,fill,minimum size=#1,inner sep=0pt,outer sep=0pt},
    dot/.default = 6pt
    }
    \node[dot] at (xP1) {};
    \node[dot] at (xP2) {};
    \node[dot] at (xP3) {};
    \node[dot] at (xP4) {};
    \node[dot] at (xP5) {};

    \draw[thick,solid] (xI23s) -- (xI13s);
    \draw[thick,solid] (xI23s) -- (xI24s);
    \draw[thick,solid] (xI24s) -- (xI45s);
    \draw[solid] (xP2) -- (xI13s); 
    \draw[solid] (xI23s) -- (xP4);
    \draw[solid] (xP2) -- (xI45s);

    \node at ($(xI24s)!0.5!(xI45s)!0.33!(xP4)$) {$T_1$};
    \node at ($(xI24s)!0.5!(xP2)!0.33!(xI45s)$) {$T_2$};
    \node at ($(xI45s)!0.5!(xP2)!0.33!(xP5)$) {$T_3$};

    \node[dot=4pt,label=right:$P_1$] at (xI45s){};
    \node[dot=4pt,label=above left:{$P_2 / P_2^0$}] at (xI24){};    
    \node[dot=4pt] at (xI23s){};
    \node[dot=4pt] at (xI13s){};

    \node[dot=2pt,label=right:$P_1^0$] at (xI45){};
    \node[dot=2pt] at (xI23){};

    \node[anchor=north] at (xP5) {$V_1^1$};
    \node[anchor=west] at (xP4) {$V_1^2 / V_2^1$};
    \node[anchor=north] at (xP2) {$V_2^2$};

    \node[anchor=west] at ( $(xP4)!0.5!(xP5)$ ) {$e_1(V_1^1,V_1^2)$};
    \node[anchor=east] at ( $(xP4)!0.45!(xP2)$ ) {$e_2(V_2^1,V_2^2)$};
    
\end{scope}

\begin{scope}[xshift=-1.5cm,scale=0.6]

 	\coordinate (xP1) at (0.0,0.0);
    \coordinate (xP2) at (1.0,0.0);
    \coordinate (xP3) at (0.5,0.86);
    \coordinate (xP4) at (1.5,0.86); 
    \coordinate (xP5) at (2.0,0.0);


    \coordinate (xI13) at ($(xP1)!0.65!(xP3)$);
    \coordinate (xI23) at ($(xP2)!0.9!(xP3)$);
    \coordinate (xI24) at ($(xP2)!0.65!(xP4)$);
    \coordinate (xI45) at ($(xP4)!0.9!(xP5)$);

    \coordinate (xI13s) at (xI13);
    \coordinate (xI23s) at ($(xP2)!0.70!(xP3)$);
    \coordinate (xI24s) at (xI24);
    \coordinate (xI45s) at ($(xP4)!0.70!(xP5)$);

    \draw[fill,color=clgr-] (xP1)--(xP2)--(xI23)--(xI13)--cycle;
    \draw[fill,color=clgr-] (xP2)--(xI24)--(xI23)--cycle;
    \draw[fill,color=clgr-] (xI45)--(xI24)--(xP2)--(xP5)--cycle;

    \draw node[rectangle,fill=clgr-,inner sep=5pt] at (0.5,1.2) {$\{\varphi_h<0\}$};

    \draw[fill,color=clgr+] (xI13) -- (xI23) -- (xP3) -- cycle;
    \draw[fill,color=clgr+] (xI23) -- (xP3) -- (xP4) -- (xI24) -- cycle;
    \draw[fill,color=clgr+] (xI24) -- (xI45) -- (xP4) -- cycle;
       
\draw node[rectangle,fill=clgr+,inner sep=5pt] at (1.5,1.2) {$\{\varphi_h>0\}$};

    \draw[very thick] (xP1) -- (xP2) -- (xP3) --cycle;
    \draw[very thick] (xP2) -- (xP4) -- (xP3) --cycle;
    \draw[very thick] (xP2) -- (xP4) -- (xP5) --cycle;

    \tikzset{
    dot/.style={circle,fill,minimum size=#1,inner sep=0pt,outer sep=0pt},
    dot/.default = 6pt
    }
    \node[dot] at (xP1) {};
    \node[dot] at (xP2) {};
    \node[dot] at (xP3) {};
    \node[dot] at (xP4) {};
    \node[dot] at (xP5) {};
\end{scope}

\end{tikzpicture}
    }
    \caption{Sketch of mesh generation using Algorithm~\ref{alg:p2s:phase2sharp}.
    On the left we sketch three triangles together with the sign of $\varphi_h$ in gray scale. 
    On the right we sketch the resulting triangulation after applying  Algorithm~\ref{alg:p2s:phase2sharp}. 
    The variables in the right triangle refer to  Algorithm~\ref{alg:p2s:phase2sharp} and describe the situation within the right triangle.
    In grey, we show the sign of $\varphi_h$, which is defined on the three large triangles indicated by the big dots. The intersection of the zero-level line and the edges is marked by  small dots.
    On both edges pointing to the top left,  these points are too close to existing vertices; based on the quality parameter $\kappa=0.3$. They are therefore shifted towards the centre of the edges, as described in line~\ref{alg:shiftLambda} of Algorithm~\ref{alg:p2s:phase2sharp}.
    The final mesh is indicated by solid lines and the newly created nodes are indicated by medium-size dots.
    }
    \label{fig:p2s:cutTriangles}
\end{figure}
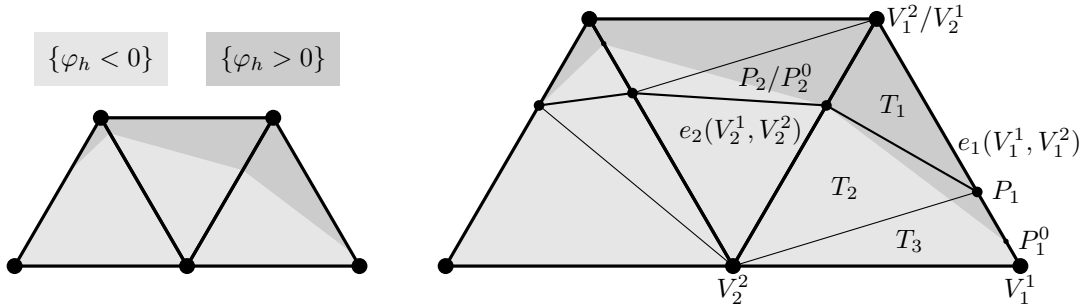

To quantify this procedure, we introduce a mesh quality parameter $\kappa \in (0,0.5]$ that quantifies how close to existing vertices newly added vertices might be created. For the precise definition see line~\ref{alg:shiftLambda} in Algorithm~\ref{alg:p2s:phase2sharp}. 
The entire process is illustrated in  Figure~\ref{fig:p2s:cutTriangles}.

Once we have generated   $\mathcal{T}_h^0$ by the described approach, we define $\Omega_0$ by those 
triangles that satisfy $\varphi_h \geq 0$ at the barycentre.

We stress that we always shift new vertices away from existing vertices (if necessary),  rather than placing them on top of existing ones.
This reduces the number of situations in which a conforming triangulation can be obtained and thus simplifies the generation of  $\mathcal T_h^0$. 
Moreover, it reduces the chance of triangles that have all three vertices on $\partial\Omega_0$; see Remark~\ref{rm:Tech:bndData}.

\begin{algorithm}
\caption{Generation of a triangulation $\mathcal T_h^0$ of $D$ that approximates $\Omega_0^\star = \{\varphi_h \geq 0\}$ }
\label{alg:p2s:phase2sharp}
\algnewcommand{\Given}{\Statex \textbf{Given: }}
\algnewcommand{\Requires}{\Statex \textbf{Requires: }}
\algnewcommand{\Returns}{\Statex \textbf{Returns: }}
\begin{algorithmic}[1]
\Given   triangulation $\mathcal T_h$, phase field $\varphi_h$ defined on $\mathcal T_h$, mesh quality parameter $\kappa \in (0,0.5]$.
\Requires $\nabla \varphi_h \not= 0$ on $\{\varphi_h=0\}$.
\Returns $\mathcal T_h^0$, i.e. triangulation of $D$ that approximates $\Omega_0^* = \{\varphi_h \geq 0\}$.
\Statex
\State Set $\mathcal T_h^0 = \emptyset$. 
\State Set $\mathcal I_h = \{ T \in \mathcal T_h\mid \exists x,y \in T\colon \varphi_h(x)\varphi_h(y) < 0\}$. 
\Comment{Triangles intersected by zero level line.}
\ForAll{$T \in \mathcal T_h$}
    \If{ $T \not \in \mathcal I_h$}
        \State $\mathcal T_h^0 := \mathcal  T_h^0 \cup T$.
        \State Continue.
    \EndIf
    
    \State Denote by $P_1^0,P_2^0$ the two intersection points of $\{\varphi_h= 0\}$ at the edges $e_1(V_1^1,V_1^2)$ and $e_2(V_2^1,V_2^2)$ of $T$ with respective vertices $V_i^j (i,j=1,2)$ of $T$.   
    \Statex\Comment{Two different points by presumption.}
    \State Calculate $\lambda_i\in (0,1)$ such that 
    $P_i^0 = \lambda_i V_i^{1} + (1-\lambda_i)V_i^{2}$, $i=1,2$, is convex combination of two vertices $V_i^{1},V_i^{2}$ of $T$. \Comment{Describe $P_i$ by scalar $\lambda_i$.}
     \State Update $\lambda_i = \max(\kappa,\min(1-\kappa,\lambda_i))$, $i=1,2$.\label{alg:shiftLambda}
         \Comment{\textit{Guarantee mesh quality.}}
     \State Set $P_i = \lambda_i V_i^{1} + (1-\lambda_i)V_i^{2}$, $i=1,2$.
     \Comment{Shift $P_i$ using $\lambda_i$.}

  \Statex
    \State Cut $T$ into triangle $T_1$ and quadrilateral $Q$ by connecting $P_1$ and $P_2$.
    \State Cut $Q$ into two triangles $T_2,T_3$ such that the smallest angle in $T_2$ and $T_3$ is maximised.
    \State $\mathcal T_h^0 := \mathcal  T_h^0 \cup \{T_1,T_2,T_3\}$.
    
\EndFor
\State \Return $\mathcal T_h^0$
\end{algorithmic}
\end{algorithm}

\begin{remark}
Algorithm~\ref{alg:p2s:phase2sharp} generates a triangulation of $D$, whose smallest angle is larger or equal to $\kappa$ times the smallest angle of the triangulation $\mathcal T_h$. Also, the change in the aspect ratio $\sigma:= \max_{T\in\mathcal T_h}\{\diam(T)\}/\min_{T\in\mathcal T_h}\{\rho(T)\}$ can be controlled in terms of $\kappa$, where $\diam(T)$ denotes the diameter of $T$ and $\rho(T)$ the radius of the inner circle of a triangle $T$.
\end{remark}

\begin{figure}
    \centering
        \includegraphics[trim = 350 260 350 60, clip,width=0.35\linewidth]{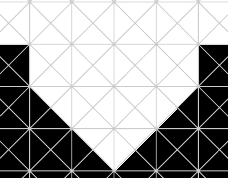}
    \hspace{2cm}        
    \includegraphics[trim = 350 260 350 60, clip,width=0.35\linewidth]{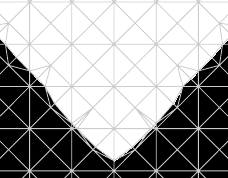}

    \caption{
        Mesh generation \mh{without (left) and with (right) inclusion of additional triangles by cutting edges, where we use $\kappa=1/8$.}  
        In both images the same phase field $\varphi_h$ is used to construct and define  $\Omega_0$ that is shown in black  together with the corresponding triangulation. On the left the values of $\varphi_h$ at barycentres of the triangles are used to assign only existing triangles of $\mathcal T_h$ to $\Omega_0$.
        On the right the approach of Algorithm~\ref{alg:p2s:phase2sharp} is used to create and define $\Omega_0$. 
        The proposed approach creates a smoother domain at the cost of locally thinner triangles.
    }
    \label{fig:tech:cutting}
\end{figure}

\begin{remark}
\label{rm:Tech:bndData}
As the suggested modification of the triangulation $\mathcal T_h$ in general leads to a reduction of the mesh quality in the sense that smaller angles and/or larger aspect ratios might appear in the triangulation $\mathcal T_h^0$, one could consider adopting the simpler approach outlined in the left-hand image of Figure~\ref{fig:tech:cutting}, which leaves the aforementioned mesh parameters unaltered during the construction of $\Omega_0$. 
In this approach $\Omega_0:=\cup_{\varphi_h(x_T) \ge 0}T$, where $x_T$ denotes the barycentre of $T$. However, this approach results in a much rougher boundary for $\Omega_0$. More importantly, some triangles may have all their vertices contained within $\partial\Omega_0$.  
In the case of Dirichlet boundary conditions and piecewise linear and continuous finite element approximations with Dirichlet boundary conditions, the state of such a triangle would be determined only by the boundary values, which could be detrimental to the approximation of the problem at hand.
\end{remark}

\subsection{Correction of geometric quantities}
\label{sec:Tech:geometricProjection}
Shape optimisation processes often also require compliance with geometric constraints.
 These are generally not satisfied by the initial value $\Omega_0$ constructed  using Algorithm~\ref{alg:p2s:phase2sharp}.
 This is due, on the one hand, to the fact that integral expressions using the phase field $\varphi_h$ are employed to describe geometric quantities (the volume of $\{\varphi_h \ge 0\}$, for example, is under-approximated by $\int_D \frac{1+\varphi_h}{2}$);
 but on the other hand also to the fact that, after applying Algorithm \ref{alg:p2s:phase2sharp}, the set $\Omega_0$ no longer coincides with $\{\varphi_h \ge 0\}$ if nodes were displaced during its execution.

In summary, we expect that geometric constraints might be violated by the proposed $\Omega_0$. 
Choosing a different level line of $\varphi_h$ in Algorithm~\ref{alg:p2s:phase2sharp} might be used to reduce this effect, however if the geometric constraint is not scalar, this may not be sufficient.

To address this, we introduce a new projection step inspired by the works \cite{MulPinRun23, HerPinSie23}, which preserve geometric quantities.
To settle notation, let us denote $\Omega_h$ to be the result of a phase field $\varphi_h$ generated by Algorithm~\ref{alg:p2s:phase2sharp}.
We seek to find $\Omega_h^*$ which is close to $\Omega_h$, but satisfies $\mathcal{G}(\Omega_h^*) = 0$, where $\mathcal{G}$ encodes the geometric constraints.
We assume that, by construction, $\mathcal{G}(\Omega_h) \approx 0$.
For some $\hat{t} \in (0,1)$, we define $\hat{V}$ to be the solution to the discrete version of 
\begin{equation}
    \hat{V} = {\rm argmin} \left\{ \frac{1}{2} \int_D  \mathcal{E} V : \mathcal{E} V \dx : V \in W^{1,\infty}_0(D;\R^d), |DV| \leq 1,\, \mathcal{G}( (\id + \hat t V)(\Omega_h)) = 0  \right\},
\end{equation}
where $\mathcal{E} V := \frac{1}{2}(DV+(DV)^T)$.
This $\hat{V}$ is numerically found as in \cite{HerPinSie23} by using a Newton method to solve a problem with the non-linear geometric constraints at each step of an iterative method which is used to solve the problem with constraint $|DV|\leq 1$.

This choice of $\hat V$ will ensure that $\Omega_h^* := ( \id + \hat{t} \hat V)(\Omega_h)$  satisfies the geometric constraint (up to a tolerance), while undergoing minimal movement.

For the subsequent shape optimisation with geometric constraints with the $W^{1,\infty}$ method, we proceed as in \cite{HerPinSie23}, minimising the action of the shape derivative, subject to the geometric constraint. This projection step only differs in the choice of functional to be minimised.


\section{Numerical investigation of exemplary PDEs}
\label{numerical-examples}
In the following we investigate topology and shape optimisation for three different physical systems, namely the Poisson equation in Section~\ref{sec:Poisson},  linear elasticity in Section~\ref{sec:Elasticity}, and finally  the Navier--Stokes equation in Section~\ref{sec:NavierStokes}.

For the sharp problems we employ the $W^{1,\infty}$ steepest descent method proposed in \cite{DecHerHin24-App}, which makes use of a hold-all domain $D$. 
We formulate all numerical examples using a hold-all domain $D$ with outer unit normal $\nu_D$ that contains the sought domain $\Omega$. 
For the phase field problems, we use the VMPT method proposed in \cite{BlankR-VMPT}.

Throughout the experiments, we set $\kappa = \frac{1}{8}$ and we move from the phase field to the sharp interface when the update is small in the sense: $\epsilon^2 \|\varphi^n - \varphi^{n-1}\|_{L^2}^2 + \epsilon\|\nabla(\varphi^n - \varphi^{n-1})\|_{L^2}^2 \leq 10^{-6}$, where $\varphi^n$ and $\varphi^{n-1}$ are consecutive iterates of the VMPT method to solve the corresponding phase field topology optimization problem.

All examples are implemented in DUNE \cite{DUNE-book}
making use of the Python bindings \cite{DunePythonReference}.


\subsection{Poisson equation}
\label{sec:Poisson}

The first example uses Poisson's equation with homogeneous Dirichlet boundary condition as PDE constraint. We use this example especially to investigate the transition from the phase field approximation to the sharp problem, outlined in Section~\ref{sec:Tech}.
Poisson's equation as constraint  has been studied in  \cite{DecHerHin25-NA}, for example, where convergence of a sequence of shapes was demonstrated under certain conditions.
Also more general objectives were considered within \cite{DecHerHin25-NA}.

We use a setup similar to \cite{DecHerHin25-NA,HerLoa23} and consider
\begin{align}
    \min J(\Omega,y) &= \int_\Omega ((y-y_d)^2 - 20 y_d) \dx
    \label{prob:sharp:poisson}
    \\
    \text{s.t. } y \in X_\Omega &\text{ satisfies} \notag
    \\
    -\Delta y &= f \mbox{ in }\Omega,\\
    y &= 0 \mbox{ on } \partial\Omega.
\end{align}
where $y_d$ and $f$ are some given functions and $X_\Omega = H_0^1(\Omega)$.

The action of the state operator $e(\Omega,y)$ applied to $p \in X_\Omega$ is given by
\begin{equation}
\label{eq:sharp:Poisson:state}
    \langle p, e(\Omega,y)\rangle
    :=
    \int_\Omega \left( \nabla y \cdot \nabla p - f p \right)\dx.
\end{equation}
\ph{We note that the existence of solutions to the shape optimization problem \eqref{prob:sharp:poisson} is quite delicate; in \cite{DZ11} a conditional existence of solutions may be provided, where the domains must satisfy a certain cone condition.}

\subsubsection{Phase field approximation}
The phase field approximation for the Poisson problem is investigated in e.g.~\cite{GarckeHKKL-2023-PhaseFieldShapeOptLaplace}.
The phase field approximation of the shape optimization problem \eqref{prob:sharp:poisson} is formulated as
\begin{align}
    \min J^\epsilon(\varphi,y^\epsilon) 
    = \int_D \frac{1+\varphi}{2} ((y^\epsilon-y_d)^2 - 20 y_d) \dx 
    &+ \gamma E^\epsilon(\varphi)\\
    \text{s.t. } (\varphi,y^\epsilon) \in H^1(D) \cap L^\infty(D) \times X_\epsilon & \text{ satisfy } &\nonumber  \\
    -\divergence( a(\varphi) \nabla y^\epsilon) + b(\varphi) y^\epsilon &= c(\varphi) f \text{ in }D,\\
    \partial_\nu y^\epsilon &= 0 \qquad \text{ on } \partial D,
\end{align}
where $\nu$ denotes the outer unit normal on $\partial D$ and $X_\epsilon = H^1(D)$.
Here we introduce interpolation functions $a,b,c \in C^{1,1}_{\text{loc}}(\R)$ defined by
\begin{align}
    \aphi(\varphi) = \cphi(\varphi) :=\ \frac{1+ \epsilon}{2} + \frac{1-\epsilon}{2}\varphi,
    \qquad
    \bphi(\varphi) :=\ \overline{b}\frac{1-\varphi}{2\epsilon^{4/3}}, 
\end{align}
that satisfy 
$a(1) = c(1) = 1$, 
$a(-1) = c(-1) = \epsilon$, 
$b(1) = 0, b(-1) = \overline{b}\epsilon^{-4/3}$.
Here $\overline{b}>0$ is a suitable scaling constant, which penalises $y^\epsilon \neq 0$ outside $\Omega$ and especially approximates the homogeneous Dirichlet data on $\partial \Omega$.
For details we refer to \cite{GarckeHKKL-2023-PhaseFieldShapeOptLaplace}.

The action of the state operator $e^\epsilon(\varphi,y^\epsilon)$ applied to $p^\epsilon \in X_\epsilon$ is given by
\begin{align}
    \langle p^\epsilon, e^\epsilon(\varphi,y^\epsilon)\rangle
    :=
    \int_D \left( \aphi(\varphi) \nabla y^\epsilon \cdot \nabla p^\epsilon + \bphi(\varphi) y^\epsilon p^\epsilon - \cphi(\varphi) f p^\epsilon \right) \dx.
\end{align}
We note that for $\gamma >0$ it follows by standard arguments, that this optimization problem admits a global solution $\varphi_* \in H^1(D) \cap L^\infty(D)$ with associated state $y_*^\epsilon \in H^1(D)$. This can be shown along the lines of \cite{GarckeHKKL-2023-PhaseFieldShapeOptLaplace}.

\subsubsection{Numerical experiment for Poisson equation}
Let us set $D = (-5/4,5/4)\times (-5/8,5/8)$ and consider $y_d(x_1,x_2) = - 0.1 (10^{-4} - (x_1^2+x_2^2)^2 + x_1^2 ) $, with $f(x) = -\Delta y_d(x)$.

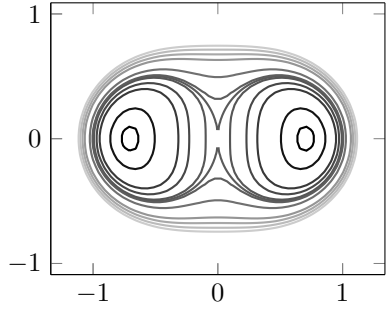
\begin{figure}
    \centering

\begin{tikzpicture}
\begin{axis}
	[width=6cm,
	domain=-1.2:1.2, 
	enlargelimits,
	view={0}{90},
	axis background/.style={fill=white},
	axis equal,
    colormap={custom}{color(0)=(black) color(1)=(black!20)}
]



\addplot3[contour prepared={labels=false},thick] file {./img/laplace/laplace-yd.tikz_contourtmp0.table};
\addplot3[contour prepared={labels=false},thick] file {./img/laplace/laplace-yd.tikz_contourtmp1.table};

\end{axis}
\end{tikzpicture}
    \caption{Isolines of the desired state $y_d$ that is used in the Poisson problem. By construction the analytical solution to \eqref{prob:sharp:poisson} is given by the levelset to one of its isolines in the form of two slightly overlapping balls.}    
    \label{fig:Poisson:yd-levelSet}
\end{figure}

This function $y_d$ has level set of approximately two slightly overlapping balls, which is visualised in Figure~\ref{fig:Poisson:yd-levelSet}.

For the phase field method we fix $\epsilon = 2^{-1}$ and $\gamma = 10^{-5}$ and we choose $\overline{b} = 100$.
The initial value is taken as $\varphi(x_1,x_2) = \frac{1}{2}( 1 + \cos(8\pi/5 x_1 ) \cos(16 \pi/5 x_2 )$.
For the spatial discretization we use a uniform triangular mesh with mesh width $h \approx 2^{-\frac{9}{2}}$.
The state spaces are discretised using globally continuous and piecewise linear finite elements.

\begin{figure}
    \centering
 
\includegraphics[trim = {25mm 71mm 143mm 10mm},clip,width=0.35\textwidth]{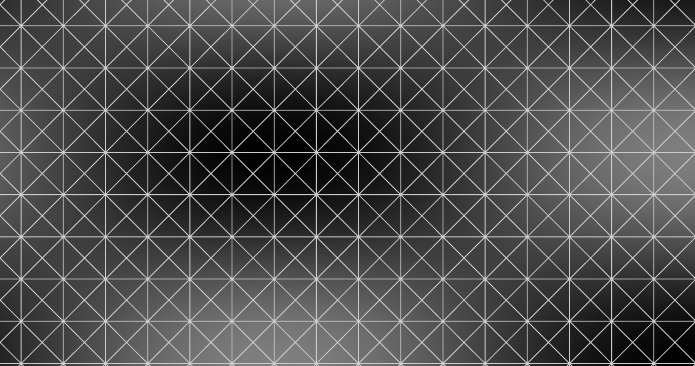}\hspace{1ex}%
\includegraphics[trim =  {143mm 71mm 25mm 10mm},clip,width=0.35\textwidth]{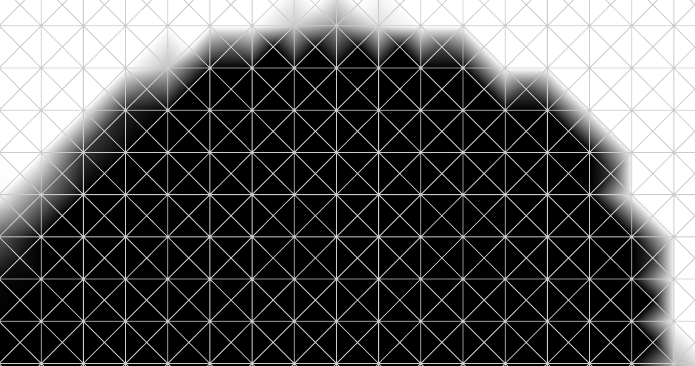}\\[1ex]%
\includegraphics[trim = {25mm 10mm 143mm 71mm},clip,width=0.35\textwidth]{img/laplace/sharp/cut/laplaceSharpinit_sharpCut_cell_mesh.pdf}\hspace{1ex}%
\includegraphics[trim =  {143mm 10mm 25mm 71mm},clip,width=0.35\textwidth]{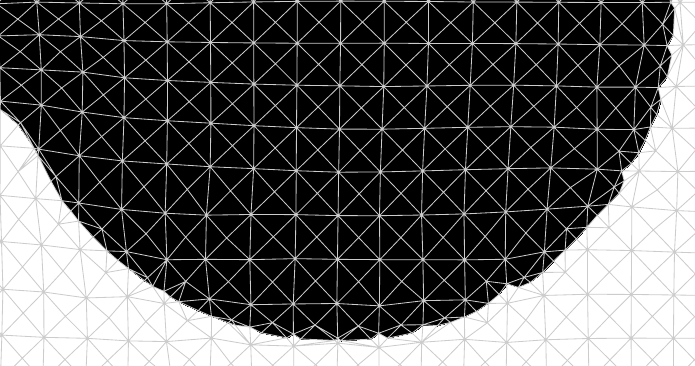}

    \caption{
    \ck{%
    Numerical solution to problem \eqref{prob:sharp:poisson}.
    On the top left we show the initial condition for the phase field problem.  Since the solution contains two slightly overlapping circles, it is symmetric with respect to the horizontal and vertical axis and we only show one quarter of every image.
    On the top left we show the final optimized phase field. On the bottom left we show the initial mesh for the sharp approach, including the additional triangles generated by the cutting procedure. On the bottom right we show the final sharp solution to the problem.
    }%
    }
    \label{fig:laplace:result-all}
\end{figure}

In Figure~\ref{fig:laplace:result-all} we show numerical results for this setup.
As the solution consists of two slightly overlapping circles, the image is symmetric with respect to the horizontal and vertical axes, so we only show one quarter of each image.

\ck{
The mesh and initial condition are shown in Figure~\ref{fig:laplace:result-all} on the top left.
After running the phase field solver until the tolerance threshold is met, we obtain the image depicted in Figure~\ref{fig:laplace:result-all} on the top right. 
We observe that the interface is poorly resolved by the considered mesh. However, this is acceptable, as the method is only used to generate an initial shape $\Omega_0$ for the second phase of the algorithm.

We then construct $\Omega_0$ as in  Algorithm~\ref{alg:p2s:phase2sharp}, the result of which is shown in Figure~\ref{fig:laplace:result-all} on the bottom left. Due to the particular structure of the optimized phase field many new cells are introduced at the boundary of $\Omega_0$.

After applying the $W^{1,\infty}$ shape optimisation method, we arrive at the domain depicted in Figure~\ref{fig:laplace:result-all} on the bottom right. In this case, most of the boundary of the shape has been flattened and is well resolved.
}

\begin{figure}
    \centering
\includegraphics[trim = {25mm 10mm 143mm 71mm},clip,width=0.35\textwidth]{img/laplace/sharp/barycenter/laplaceSharpinit_barycenter_cell_mesh.pdf}\hspace{1ex}%
\includegraphics[trim =  {143mm 10mm 25mm 71mm},clip,width=0.35\textwidth]{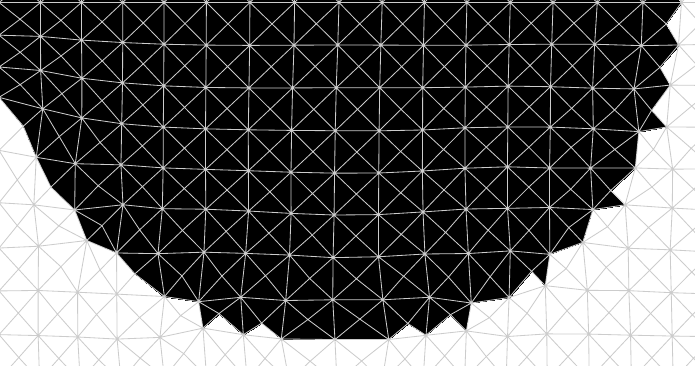}
    \caption{On the left we show $\Omega_0$  generated from the phase field in Figure~\ref{fig:laplace:result-all} (top right) and without cutting procedure and  on the right the corresponding optimized shape. We observe that the shape optimization algorithm is not able to fully remove the spikes on the boundary of the shape.}
    \label{fig:laplace:sharp:noCut}
\end{figure}

\ck{%
For comparison we show in Figure~\ref{fig:laplace:sharp:noCut} the numerical result for the case that the cutting procedure from Algorithm~\ref{alg:p2s:phase2sharp} is not used and instead values are assigned via the barycentre value. On the left we show $\Omega_0$ and on the right the optimized shape, starting with this $\Omega_0$.
As expected, we observe that the poor resolution of the boundary of $\Omega_0$ is mainly preserved throughout the optimization, especially compared to the result with cutting in Figure~\ref{fig:laplace:result-all} (bottom right), and is only poorly flattened. 
}%

\begin{figure}
    \centering
\includegraphics[trim = {25mm 10mm 143mm 71mm},clip,width=0.35\textwidth]{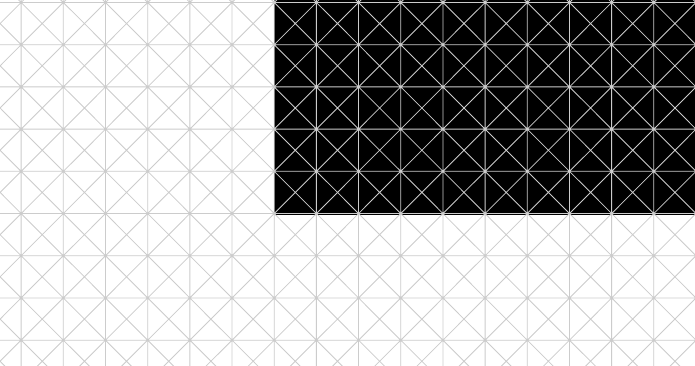}\hspace{1ex}%
\includegraphics[trim =  {143mm 10mm 25mm 71mm},clip,width=0.35\textwidth]{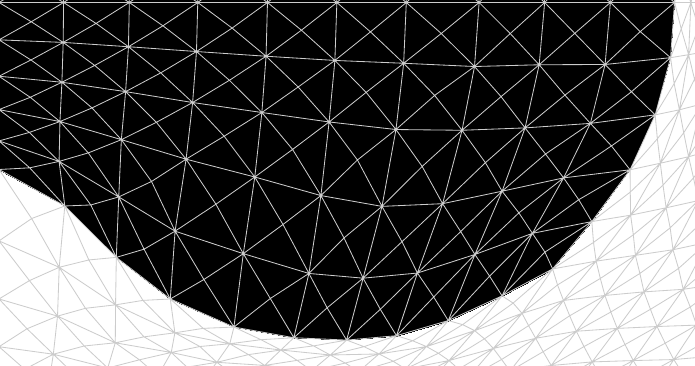}
    \caption{%
    \ck{%
    The numerical results using an initial guess $\Omega_0$ for the sharp method, that is significantly different from the optimized shape.
    On the left we show the initial guess $\Omega_0$ and on the right the resulting optimized shape. Since $\Omega_0$ is a poor initial guess, we observe strong mesh deterioration. As can be seen in Figure~\ref{fig:laplace:result-all} (bottom right) this can be prevented by using a more suitable initial guess. 
    }%
    }
    \label{fig:laplace:guess}
\end{figure}

In Figure \ref{fig:laplace:guess}, we depict the result of an initial guess that $\Omega_0 = (-5/8,5/8) \times (-5/16,5/16)$. One might consider this a poor initial guess.
We observe that the well resolve boundary of $\Omega_0$ is preserved, while the mesh is significantly distorted due to the initial guess being significantly different from the optimised shape.

We note that when we used the phase field method, the shape was already roughly correct, so the shape optimisation procedure did not have much to do.
\ck{%
On the other hand, we observe that the coarse initial mesh generates a phase field that leads to a rather rough initial shape $\Omega_0$ that is not fully smoothed during the optimization process.
}%



\subsection{Elasticity problem}
\label{sec:Elasticity}
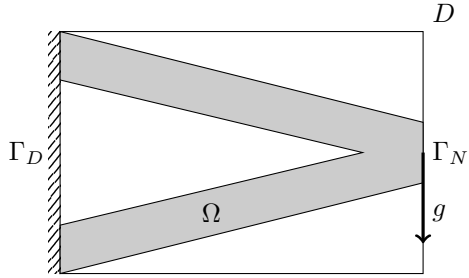
\begin{figure}
    \centering
    \begin{tikzpicture}[scale=0.8]
    
    \begin{scope}
    \pgfmathsetmacro{\LL}{6.0}
    \pgfmathsetmacro{\HH}{4.0}
    
        \draw (0.0,0.0) rectangle (\LL,\HH);
  
        \draw[fill=black!20](0.0,0.0) -- (6.0,1.5) -- (6.0,2.5) -- (0.0,4.0) --  (0.0,3.2) -- (5.0,2.0) -- (0.0,0.8) -- cycle;
        
        \draw (2.5,1.0) node {$\Omega$};
        
        \draw[very thick,->] (\LL,2.0) -- (\LL,0.5);
        
\fill[pattern={Lines[angle=50, distance=2.5pt,line width=0.5pt]}, pattern color=black!100] (-0.2,0.0) rectangle (0.0,\HH);

         \draw (-0.1,2) node[anchor=east] {$\Gamma_D$};
        \draw (6,2.0) node[anchor=west] {$\Gamma_N$};

        \draw (6,1.0) node[anchor=west] {$g$};

        \draw (6,4) node[anchor=south west] {$D$};
        
        \end{scope}

\end{tikzpicture}
    \caption{Setup for the elastic example. The unknown structure is fixed on the left at $\Gamma_D$ and  a force $g$ is applied on the right at $\Gamma_N$. The gray domain $\Omega$ is shown for illustration purposes.}
    \label{fig:elastic:setup}
\end{figure}

Perhaps the most interesting application of our strategy is the task of minimising the elastic compliance. 
It is well known that elastic compliance has non-trivial topology; as such, standard sharp-shape optimisation methods may fail since topology changes may occur during the optimisation process, as indicated in Figure ~\ref{fig:exampleTopOpt} of the introduction.
We consider a situation similar to the Messerschmitt--Bölkow--Blohm beam example from \cite[\S~4.6]{PapFarSur-2021-MultipleSolutionsTopoOpt}, i.e. we consider a structure $\Omega$, that is fixed on the left of $D$ at $\Gamma_D$, while on the right at $\Gamma_N$ a force $g$ pushes down from above; see Figure~\ref{fig:elastic:setup} for a sketch of the situation, including a possible structure $\Omega$.

The problem we consider may be written as: Given a hold-all domain $D$ with outer unit normal $\nu_D$ and  with $\Gamma_N \subset \partial D$ and $\Gamma_D \subset \partial D$. 
Seek $\Omega \subseteq D$ with outer unit normal $\nu_\Omega$, $\Gamma_N \subset \partial\Omega$ and $\Gamma_D \cap \partial \Omega \neq \emptyset$ such that
\begin{align}
    \label{prob:sharp:elastic}
    \min J(\Omega,y) &= \int_{\Gamma_N} g\cdot y \ds  
    + \gamma|\partial \Omega \setminus \partial D| \\
    \text{s.t. } y \in X_\Omega& \text{ satisfies } \notag    \\
    -\divergence{\sigma} &= 0 \mbox{ in } \Omega,\\
    \sigma& = 2 \mu \mathcal E(y) + \lambda \trace(\mathcal E (y))I,\\
    \sigma \nu_\Omega &= g \mbox{ on } \Gamma_N,\\
    y &= 0 \mbox{ on } \Gamma_D,\\
    \sigma \nu_\Omega &= 0 \mbox{ on } \partial \Omega \setminus (\Gamma_N\cup \Gamma_D),
    \\
    |\Omega| &= c,
\end{align}
where $\mathcal E(y) := \frac{1}{2} (D y + (D y)^T)$ is the symmetric gradient, $\mu>0$ is the shear parameter, $\lambda$ is the first Lam\'e parameter, $c$ is a chosen volume constraint and $X_\Omega:= \{ v \in H^1(\Omega;\R^d) : v|_{\Gamma_D} = 0\}$.
The action of the state operator $e(\Omega,y)$ applied to $p \in X_\Omega$ is given by
\begin{equation}
\label{eq:sharp:elastic:state}
    \langle p, e(\Omega ,y)\rangle
    :=
    \int_\Omega \left(
    2 \mu \mathcal{E}(y) : \mathcal{E} (p)
    +
    \lambda \divergence{y} \divergence{p}
    \right) \dx
    -
    \int_{\Gamma_{N}} g \cdot p\ \ds.
\end{equation}

\subsubsection{Phase field regularisation}
The phase field regularised problem is discussed e.g. in \cite{GarckeHuettlKnopf-2022-TopoOptEigValElasticity,GarckeHKK-2024-MultiPhaseElasticSpectral} and is given by
\begin{align}
    \min J^\epsilon(\varphi,y^\epsilon) = \int_{\Gamma_N} g\cdot y^\epsilon \ \ds 
    &+ \gamma E^\epsilon(\varphi)\\
    \text{s.t. } (\varphi, y^\epsilon)  \in  H^1(\Omega)\cap L^\infty(\Omega) \times  X_\epsilon &\text{ satisfy }\\
        -\divergence{\sigma} &= 0 \mbox{ in } D,\\
    \sigma& =\alpha(\varphi) \left(2 \mu \mathcal E(y^\epsilon) + \lambda \trace(\mathcal E (y^\epsilon))I \right),\\
    \sigma \nu_D &= g \mbox{ on } \Gamma_N,\\
    y^\epsilon &= 0 \mbox{ on } \Gamma_D,\\
    \sigma \nu_D &= 0 \mbox{ on } \partial D \setminus (\Gamma_N\cup \Gamma_D),
    \\
    \int_D \frac{1+\varphi}{2} dx &= c,
\end{align}
where $\mathcal{E}$ is again the symmetric gradient,  $\mu$, $\lambda$ are the Lam\'e parameters and $X_\epsilon = \{ v \in H^1(D;\R^d) : v|_{\Gamma_D} = 0\}$.
Moreover for $\underline{\alpha} >0$ 
\begin{align}
    \alpha(\varphi) = 
    \underline \alpha \epsilon 
    + (1-\underline \alpha \epsilon) \left(\frac{1+\varphi}{2}\right)^2
\end{align}
 models a very weak material outside of the domain $\Omega$ by scaling the material properties of $\Omega$ with $\underline \alpha \epsilon$.
The action of the state operator $e^\epsilon(\varphi,y^\epsilon)$ applied to $p^\epsilon \in X_\epsilon$ is given by
\begin{equation}
\label{eq:diffuse:elastic:state}
    \langle p^\epsilon, e(\varphi,y^\epsilon)\rangle
    :=
    \int_D \aphi(\varphi) \left( \mu \mathcal{E} y^\epsilon : \mathcal{E} p^\epsilon  
    + \lambda \divergence{y^\epsilon} \divergence{p^\epsilon} \right) \dx 
    - \int_{\Gamma_{N}} g \cdot p^\epsilon\ \ds.
\end{equation}
Existence of optimal controls and necessary first order optimality conditions for this problem are, e.g., shown in \cite{BlankEtAl-2012-PhaseFieldStructuralTopoOpt}, including sharp interface asymptotics for $\epsilon\to 0$.

\subsubsection{Numerical experiments for elasticity}
 We consider a setup, similar to \cite[\S~4.6]{PapFarSur-2021-MultipleSolutionsTopoOpt}.
 The parameters are 
 $D = (-1,1)\times(-0.5,0.5)$, 
 $\Gamma_N = \{x_1 = 1, x_2 \in(-0.125,0.125)\}$, 
 $\Gamma_{D} = \{x_1=-1\}$, $g =(0,10)^T $, 
 and for the material the Lam\'e parameters are $\mu = 75.38$, $\lambda = 64.62$.
For the phase field problem, we set $\epsilon = 2^{-5}$, $\gamma = 2^{-5}$, and $\underline \alpha=0.05$, 

Moreover, in the phase field setting, we introduce the additional constraint that $\varphi\geq 0$ on $\Gamma_N$.
This ensures that no relevant boundary is missed when transferring from the diffuse method to the sharp method.

The phase field approach inherently introduces surface area regularization with  strength $\gamma$.
When performing shape optimisation, we split into a well-posed surface area regularised case and an ill-posed case without surface area regularisation.
This allows us to determine whether pointed shapes are obtained as a result of the sharp setting.

\begin{figure}
    \centering
    \begin{subfigure}[t]{.45\linewidth}
        \centering
        \includegraphics[width = \linewidth]{img/elastic_reg/phase/init/elasticPhaseinit_cell_mesh.pdf}
        \caption{Initial $\varphi$ for the diffuse interface method.}
        \label{fig:elastic:diffuse:initial}
    \end{subfigure}
    \hfill
    \begin{subfigure}[t]{.45\linewidth}
        \centering
        \includegraphics[width = \linewidth]{img/elastic_reg/phase/final/elasticPhasefinal_cell_mesh.pdf}
        \caption{Final $\varphi$ for the diffuse interface method.}
        \label{fig:elastic:diffuse:final}
    \end{subfigure}
    \\
    \begin{subfigure}[t]{.45\linewidth}
        \centering
        \includegraphics[width = \linewidth]{img/elastic_reg/sharp/init/elasticSharpinit_cell_mesh.pdf}
        \caption{The cut mesh (with incorrect volume).}
        \label{fig:elastic:sharp:preProjection}
    \end{subfigure}
    \hfill
    \begin{subfigure}[t]{.45\linewidth}
        \centering
        \includegraphics[width = \linewidth]{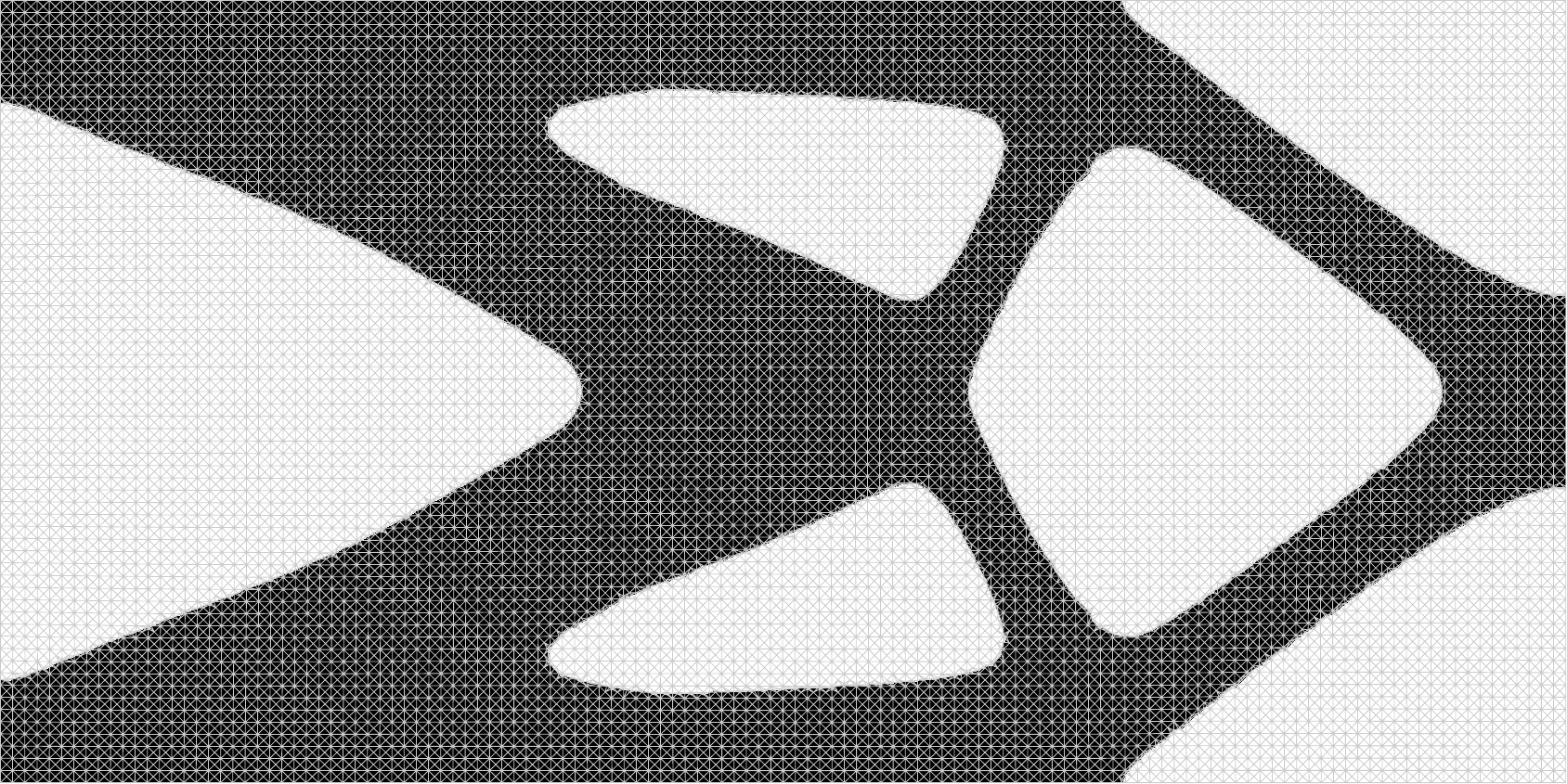}
        \caption{Initial domain for the sharp method. (could possibly overlay the cut mesh on top)}
        \label{fig:elastic:sharp:initial}
    \end{subfigure}
    \\
    \begin{subfigure}[t]{.45\linewidth}
        \centering
        \includegraphics[width = \linewidth]{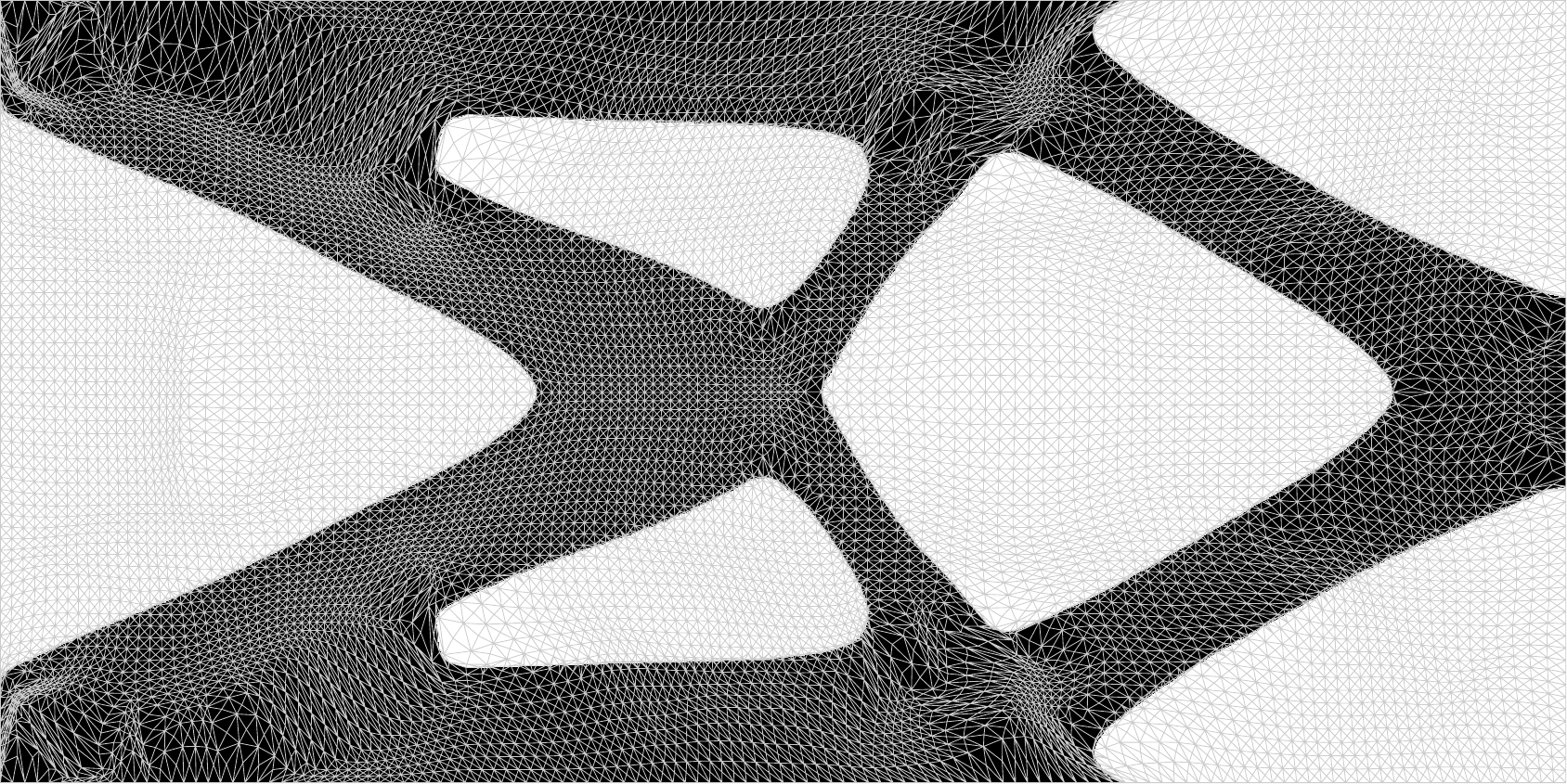}
        \caption{Final domain for the sharp method with surface tension.}
        \label{fig:elastic:sharp:finalWithTen}
    \end{subfigure}
    \hfill
    \begin{subfigure}[t]{.45\linewidth}
        \centering
        \includegraphics[width = \linewidth]{img/elastic_noReg/sharp/final/elasticSharpfinal_cell_mesh.pdf}
        \caption{Final domain for the sharp method without surface tension.}
        \label{fig:elastic:sharp:finalNoTen}
    \end{subfigure}
    \caption{Numerical solution to problem \eqref{prob:sharp:elastic}.
            At the top we show the initial and optimised phase function.
            In the next row we show the result of Algorithm \ref{alg:p2s:phase2sharp} and the projection of this, which is the initial condition for the sharp setting.
            In the final row, we show the result of the $W^{1,\infty}$ shape optimisation with and without surface tension. We observe that the symmetric voids in the middle are more pronounced in the case without surface tension.}
    \label{fig:experiment:elastic}
\end{figure}

\ck{
The results of our numerical test are shown in Figure~\ref{fig:experiment:elastic}.
We start the optimization with the constant initial condition of $\varphi = 0$, this is illustrated in Figure~\ref{fig:elastic:diffuse:initial}.
The optimized phase field is shown in Figure~\ref{fig:elastic:diffuse:final}.
After applying Algorithm~\ref{alg:p2s:phase2sharp} we obtain a first initial shape for the sharp method, which is shown in Figure~\ref{fig:elastic:sharp:preProjection}. This initial shape generally does not satisfy the geometric constraint.
We then follow the procedure outlined in Section~\ref{sec:Tech:geometricProjection} which leads to minor changes and the actual initial shape. This is depicted in Figure~\ref{fig:elastic:sharp:initial} for the sharp method.

At this point, we run two shape optimisations using $W^{1,\infty}$ method, firstly with positive surface area penalty, $\gamma = 2^{-5}$, which results in the optimized shape shown in Figure~\ref{fig:elastic:sharp:finalWithTen}, and with no penalty, $\gamma = 0$, which yields the optimized shape shown in Figure~\ref{fig:elastic:sharp:finalNoTen}.
For the case with area regularisation, we observe that the two holes in the middle of the optimised shape in Figure~\ref{fig:elastic:sharp:finalWithTen} are very close to the optimized phase field solution in Figure~\ref{fig:elastic:diffuse:final}, while the main hole on the left is more stretched. 
The same applies to the optimisation without area regularisation, in which the two holes are also more pronounced. This indicates that the thinner, longer holes have only a minor effect on the objective, which is surpassed by the area regularisation in Figure~\ref{fig:elastic:sharp:finalWithTen}.
}

\ck{
\begin{remark}
    In this numerical example, we obtain similar results regardless of whether we use surface area regularisation. 
   However, when using the phase field approach to generate an initial shape, surface area regularisation is introduced to the problem in general.
   Therefore, ignoring this in the $W^{1,\infty}$ method, might lead to the phase field generating an incorrect initial topology for the problem.
\end{remark}
}

\subsection{Navier--Stokes equation}
\label{sec:NavierStokes}
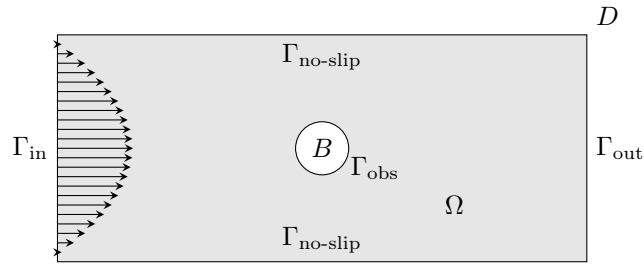
\begin{figure}
    \centering
    {
    \tikzset{every picture/.style={scale=0.5}}
    \begin{tikzpicture}[scale=1.0]

    \pgfmathsetmacro{\LL}{7.0}
    \pgfmathsetmacro{\HH}{3.0} 
    
		\draw[fill=black!10] (-\LL,-\HH) rectangle (\LL,\HH);

        \node[circle,draw,fill=white,minimum size=0.5] (B) at (0,0) {$B$};	
        \draw (B.south east) node[anchor=west] {$\Gamma_{\text{obs}}$};

        \draw (\LL,\HH) node[anchor=south west] {$D$};
        
		\draw (-\LL,0.0) node[anchor=east] {$\Gamma_{\text{in}}$};
		\draw (\LL,0.0) node[anchor=west] {$\Gamma_{\text{out}}$};
		
		\draw (0.0,\HH) node[anchor=north] {$\Gamma_{\text{no-slip}}$};
		\draw (0.0,-\HH) node[anchor=south] {$\Gamma_{\text{no-slip}}$};
		\draw (\LL/2,-\HH/2) node {$\Omega$};
		\foreach \y in {-2.75,-2.5,...,2.75}{
				\draw [thin,-stealth] (-\LL,\y) -- (-5-0.25*\y*\y,\y);
			};
\end{tikzpicture}

    
    }
    \caption{Geometric setup of the obstacle in Navier--Stokes fluid example. Fluid enters on the left side of the domain $D$, flows around an obstacle $B$ and leaves the domain on the right side of the domain.}
    \label{fig:fluid:setup}  
\end{figure}

The problem we consider here is a classical example of fluid-mechanical shape optimization and is considered in e.g.~\cite{MulKuhSie21}.
We consider an obstacle $B$  in a Navier--Stokes flow, where the fluid enters the domain on the left through the inlet boundary $\Gamma_{\text{in}}$ and leaves the domain on the right  through the outlet boundary $\Gamma_{\text{out}}$. See Figure~\ref{fig:fluid:setup} for a sketch.

The objective of the problem is to minimise the energy dissipation due to viscosity around an obstacle with a fixed centre of mass and volume.
In this setting, we note that fixing the centre of mass of the fluid domain $\Omega$ directly fixes the centre of mass of the obstacle.

The problem is stated as:
Seek $\Omega$ such that $\partial D \subset \partial \Omega$,
\begin{align}
    \min J(\Omega,y) = \frac{\mu}{2} \int_\Omega Dy : Dy\ \dx \label{prob:NS:sharp} \\
    \text{s.t. } (y,\eta) \in X_\Omega &\text{ satisfy }\notag \\
    -\mu\Delta y + (y\cdot\nabla) y + \nabla \eta &= 0 \text{ in }\Omega,\\
    -\divergence{y} &= 0 \mbox{ in } \Omega,\\
    y &= g \text{ on } \Gamma_{\text{in}},\\
    y &= 0 \text{ on } \Gamma_{\text{no-slip}},\\
    ( \mu (\nabla y) ^T+ pI ) \nu_D &= 0 \text{ on } \Gamma_{\text{out}},\\
    |\Omega| = c, \int_\Omega x\ \dx &= 0. \label{eq:NS:geoconstr}
\end{align}
Here $\mu>0$ denotes the viscosity of the fluid, $g$ suitable Dirichlet boundary data  and  the state space $X_\Omega$ is given by 
\begin{align*}
    X_\Omega:= \{ (y,\eta) \in H^1(\Omega;\R^d) \times L^2(\Omega) : 
    y|_{\Gamma_{\text{in}}} = g,\,
    y|_{\Gamma_{\text{no-slip}}} = 0,\, 
    \int_\Omega \eta\ \dx = 0 \}.
\end{align*}
The geometric constraints \eqref{eq:NS:geoconstr} fix the volume of the fluid domain $\Omega$ to a given value $c$ and fix the centre of mass of $\Omega$.
The action of the state operator $e(\Omega,(y,\eta))$ applied to $(p,\zeta) \in X_\Omega^0$ is given by
\begin{equation}
\label{eq:sharp:fluid:state}
    \langle (p,\zeta), e(\Omega,(y,\eta))\rangle
    :=
    \int_\Omega \left( \nu D y : D p + \left( y \cdot \nabla \right) y p + \eta \divergence{p} \right) \dx
    +
    \int_\Omega \zeta \divergence{y} \ \dx.
\end{equation}
where
\begin{align*}
    X_\Omega^0:= \{ (y,\eta) \in H^1(\Omega;\R^d) \times L^2(\Omega) : 
    y|_{\Gamma_{\text{in}}} = 0,\,
    y|_{\Gamma_{\text{no-slip}}} = 0\}.
\end{align*}

\subsubsection{Phase field regularisation}
The phase field regularised problem is e.g.~discussed in 
\cite{GarHHk-2015-NumApproxPhaseTopoOptFluids,GarHKL-2018-TopOptPF-StateConstraints} and
is given as
\begin{align}
    \min J^\epsilon(\varphi,y^\epsilon) = 
    \frac{\mu}{2} \int_D D y^\epsilon : Dy^\epsilon \ \dx
    + \frac{1}{2}\int_D\alpha(\varphi) | y^\epsilon|^2 \dx &
    +\gamma E^\epsilon(\varphi) \\
    \text{s.t. }(\varphi,(y^\epsilon,\eta^\epsilon)) \in H^1(D) \cap L^\infty(D) \times X_\epsilon & \text{ satisfy }\notag\\
    \alpha(\varphi) y^\epsilon  -\mu\Delta y^\epsilon + (y^\epsilon\cdot\nabla) y^\epsilon + \nabla \eta^\epsilon &= 0\text{ in } D ,\\
    -\divergence{y^\epsilon} &= 0 \text{ in } D,\\
    y^\epsilon &= g \text{ on }  \Gamma_{\text{in}},\\
    y^\epsilon &= 0 \text{ on }  \Gamma_{\text{no-slip}},\\
    ( \mu (\nabla y^\epsilon) ^T+ \eta^\epsilon I ) \nu_D &= 0 \text{ on } \Gamma_{\text{out}},\\
    \int_D \frac{1+\varphi}{2} \dx &= c, \label{eq:NS-PF:vol}\\ 
    \int_D \frac{1+\varphi}{2} x\ \dx &= 0. \label{eq:NS-PF:com}
\end{align}
where again $\mu$ denotes the viscosity of the fluid and $g$ given Dirichlet boundary data. 
Moreover $\alpha(\varphi) = \overline{\alpha}\epsilon^{-1}(1-\varphi)$ with a fixed positive constant $\overline \alpha$, such that $\alpha(1) = 0$ and $\alpha(-1) = 2\overline \alpha \epsilon^{-1}$ penalizes $y \neq 0$ inside the obstacle; comparable to the function $b$ in the Poisson example.
The state space is given by 
\begin{align*}
    X_\epsilon:= \{ (y^\epsilon,\eta^\epsilon) \in H^1(D;\R^d) \times L^2(D) : 
    y|_{\Gamma_{\text{in}}} = g,\,
    y|_{\Gamma_{\text{no-slip}}} = 0,\, 
    \int_D \eta\ \dx = 0 \}.
\end{align*}
The action of the state operator $e^\epsilon(\varphi,(y^\epsilon,\eta^\epsilon))$ applied to $(p^\epsilon,\zeta^\epsilon) \in X_\epsilon^0$ is given by
\begin{align*}
    \langle (p^\epsilon,\zeta^\epsilon) , e(\varphi,(y^\epsilon,\eta^\epsilon) )\rangle
    :=
    \int_D \left( \mu D y^\epsilon : D p^\epsilon + (y^\epsilon \cdot \nabla) y^\epsilon \cdot p^\epsilon + \alpha(\varphi) y^\epsilon p^\epsilon - \eta^\epsilon \divergence p^\epsilon + \zeta^\epsilon \divergence y^\epsilon\right) \dx,
    \label{eq:NS-PF:stateweak}
\end{align*}
where 
\begin{align*}
    X_\epsilon^0:= \{ (p^\epsilon,\zeta^\epsilon) \in H^1(D;\R^d) \times L^2(D) : 
    y|_{\Gamma_{\text{in}}} = 0,\,
    y|_{\Gamma_{\text{no-slip}}} = 0,\, 
    \int_D \eta\ \dx = 0 \}.
\end{align*}
Concerning existence of optimal solutions to this optimization problem  we refer to \cite{GarHHk-2015-NumApproxPhaseTopoOptFluids}.

\subsubsection{Numerical experiments for Navier--Stokes}
We use the computational domain $D=(-5,5)\times(-2.5,2.5)$ with inflow $g(x_1,x_2) = (1-(x_2/2.5)^2,0)$ on $\Gamma_{\text{in}} = \{x_1=-5\}$, and open boundary outflow at $\Gamma_{\text{out}} = \{x_1=5\}$
The viscosity is $\mu=1$.
Moreover, we fix the centre of mass at $(0,0)$ and the volume of the obstacle to $c=2.25$.
Additionally, for the phase field approximation we use $\overline \alpha=50$, $ \epsilon = 0.5$, $\gamma=0.5$.
As initial phase field we use $\varphi(x_1,x_2) = 0.91 + 10^{-2} \cos\left(\frac{\pi x_1}{5}  \right) + 7 \times 10^{-2} \cos\left(\frac{2\pi x_2}{5} \right) $.
Fixing the volume of the obstacle fixes the mean value of $\varphi$ to be $\frac{1}{|D|}\int_D\varphi\,\dx = 0.91$.
To avoid phase accumulating at the boundary, we restrict to $\varphi|_{\partial D} = 1$.

\begin{figure}
    \centering
    \begin{subfigure}[t]{.45\linewidth}
        \centering
        \includegraphics[trim=180 100 180 100,clip,width = \linewidth]{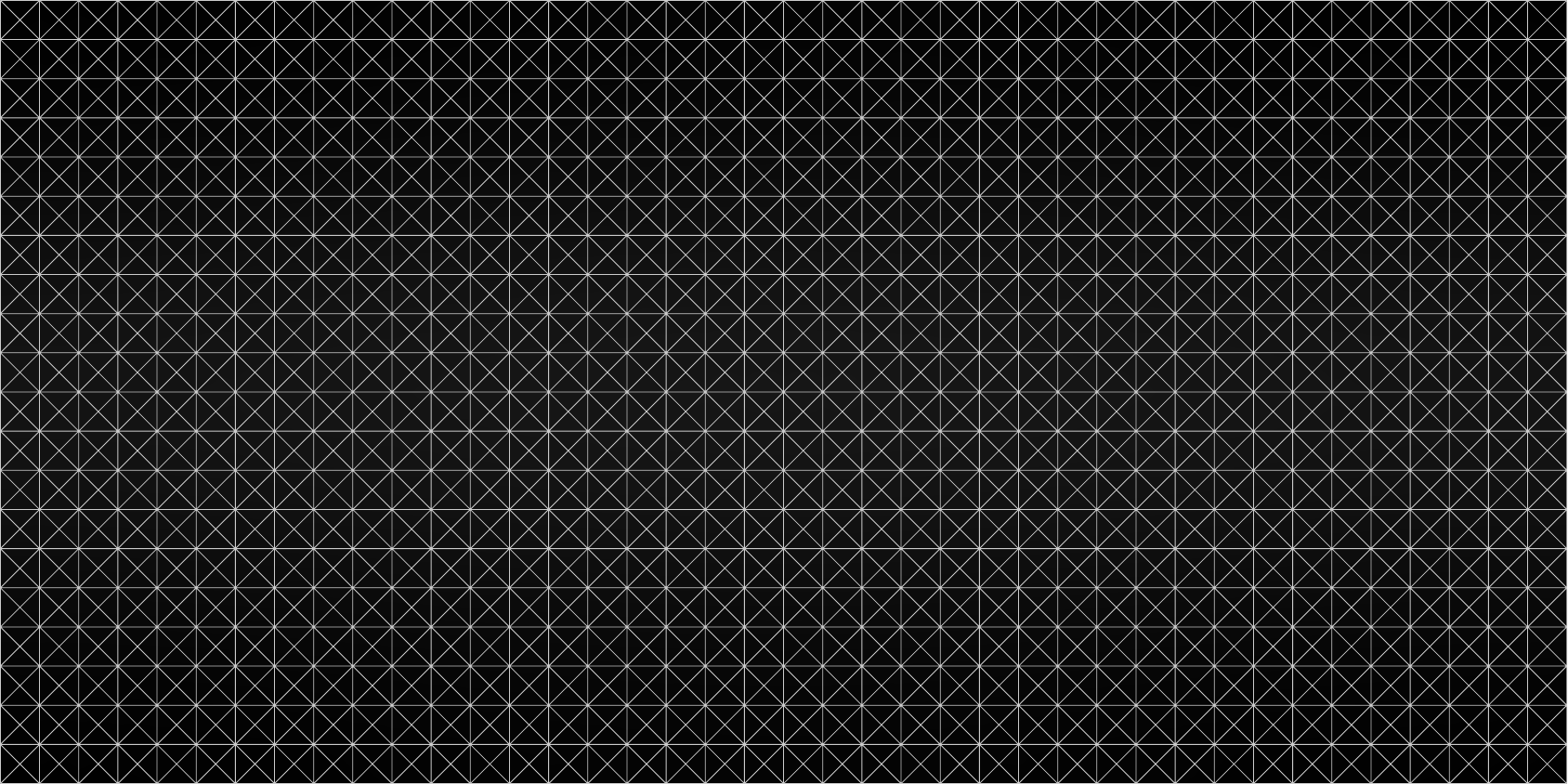}
        \caption{Initial $\varphi$ for the diffuse interface method (blow up).}
        \label{fig:navier_stokes:diffuse:initial}
    \end{subfigure}
    \hfill
    \begin{subfigure}[t]{.45\linewidth}
        \centering
        \includegraphics[trim=180 100 180 100,clip,width = \linewidth]{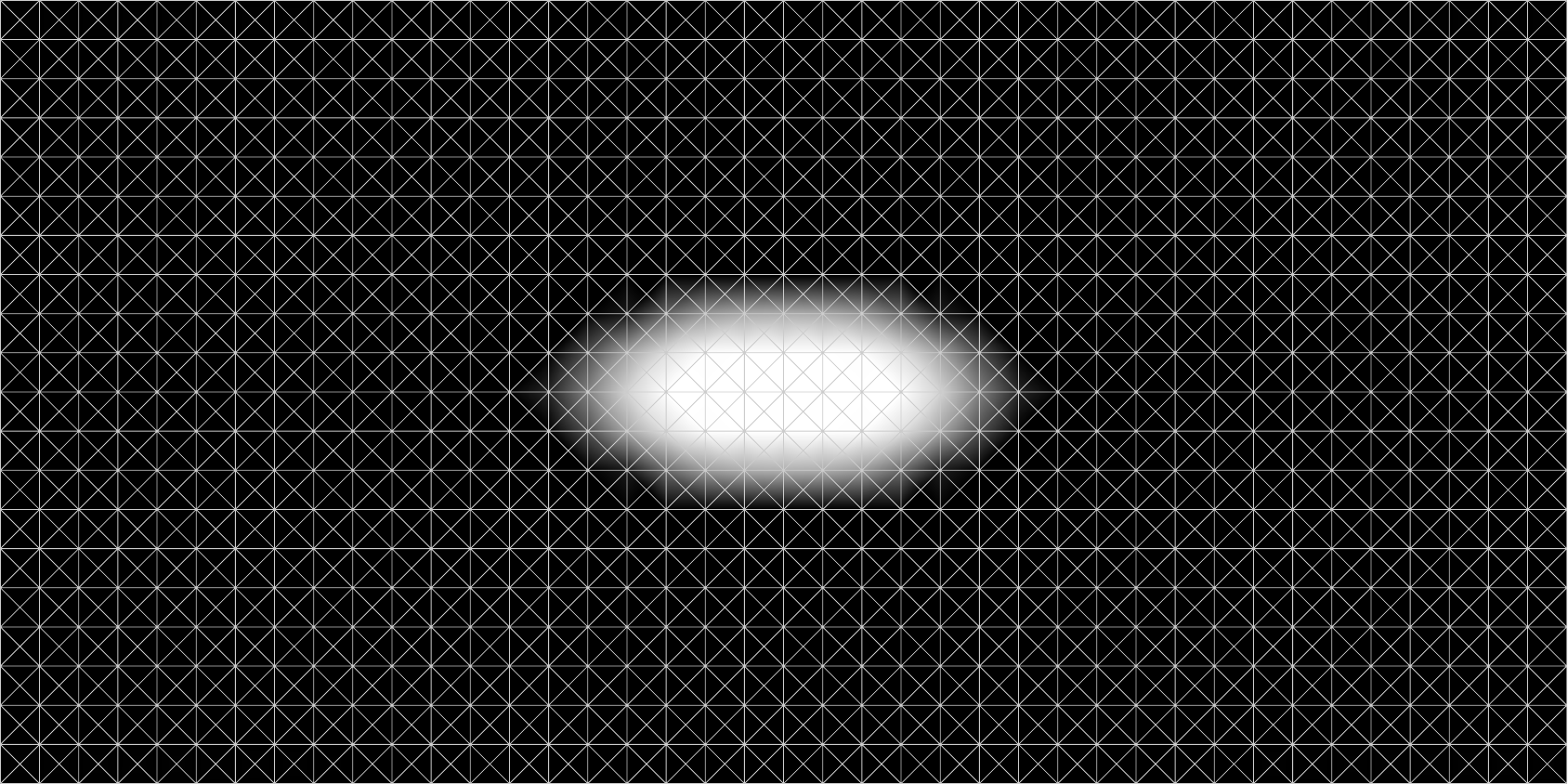}
        \caption{Final $\varphi$ for the diffuse interface method (blow up).}
        \label{fig:navier_stokes:diffuse:final}
    \end{subfigure}
    \\
    \begin{subfigure}[t]{.45\linewidth}
        \centering
        \includegraphics[trim=180 100 180 100,clip,width = \linewidth]{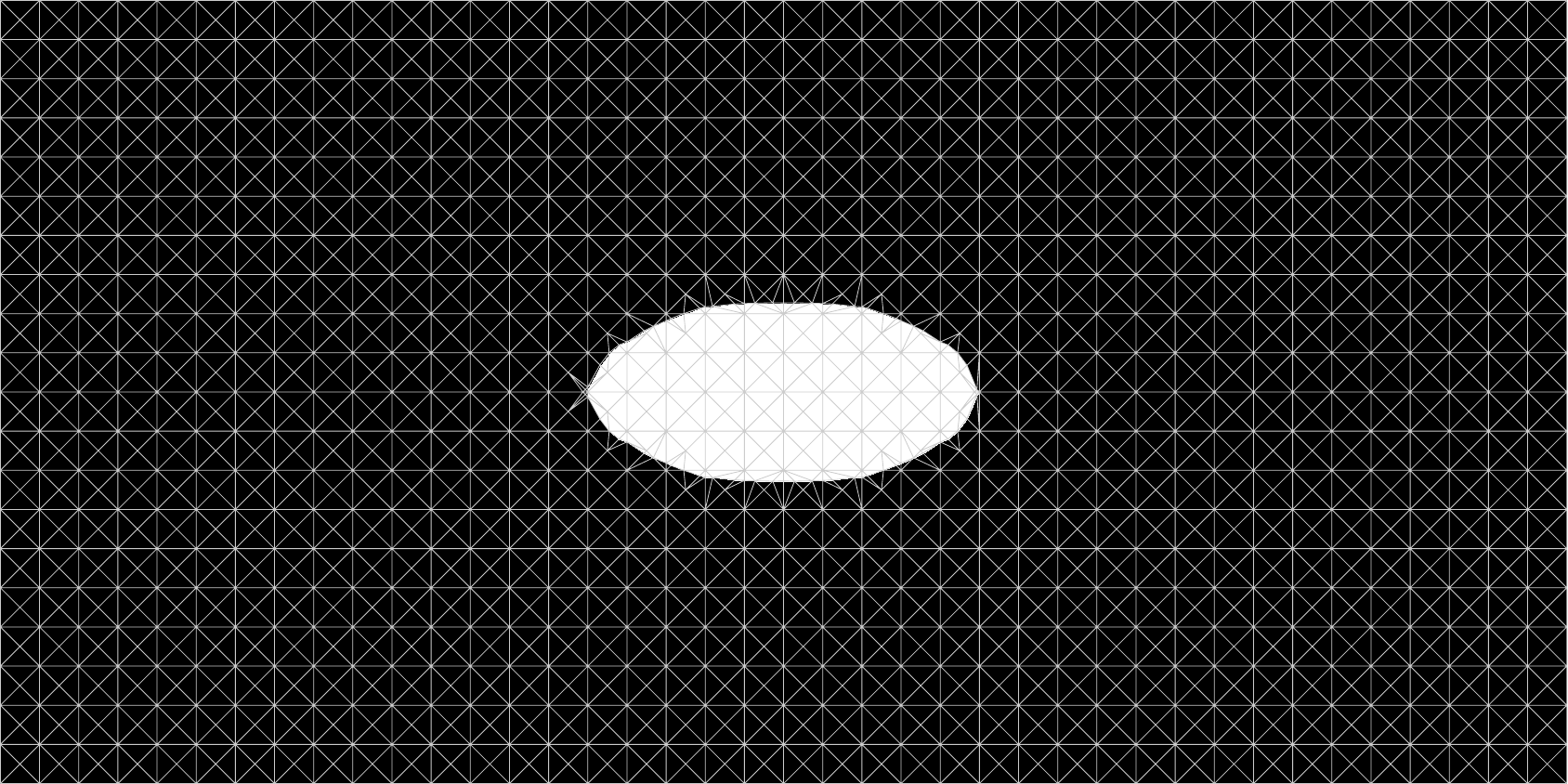}
        \caption{The cut mesh; with incorrect volume and barycentre (blow up).}
        \label{fig:navier_stokes:sharp:initial}
    \end{subfigure}
    \hfill
    \begin{subfigure}[t]{.45\linewidth}
        \centering
        \includegraphics[trim=180 100 180 100,clip,width = \linewidth]{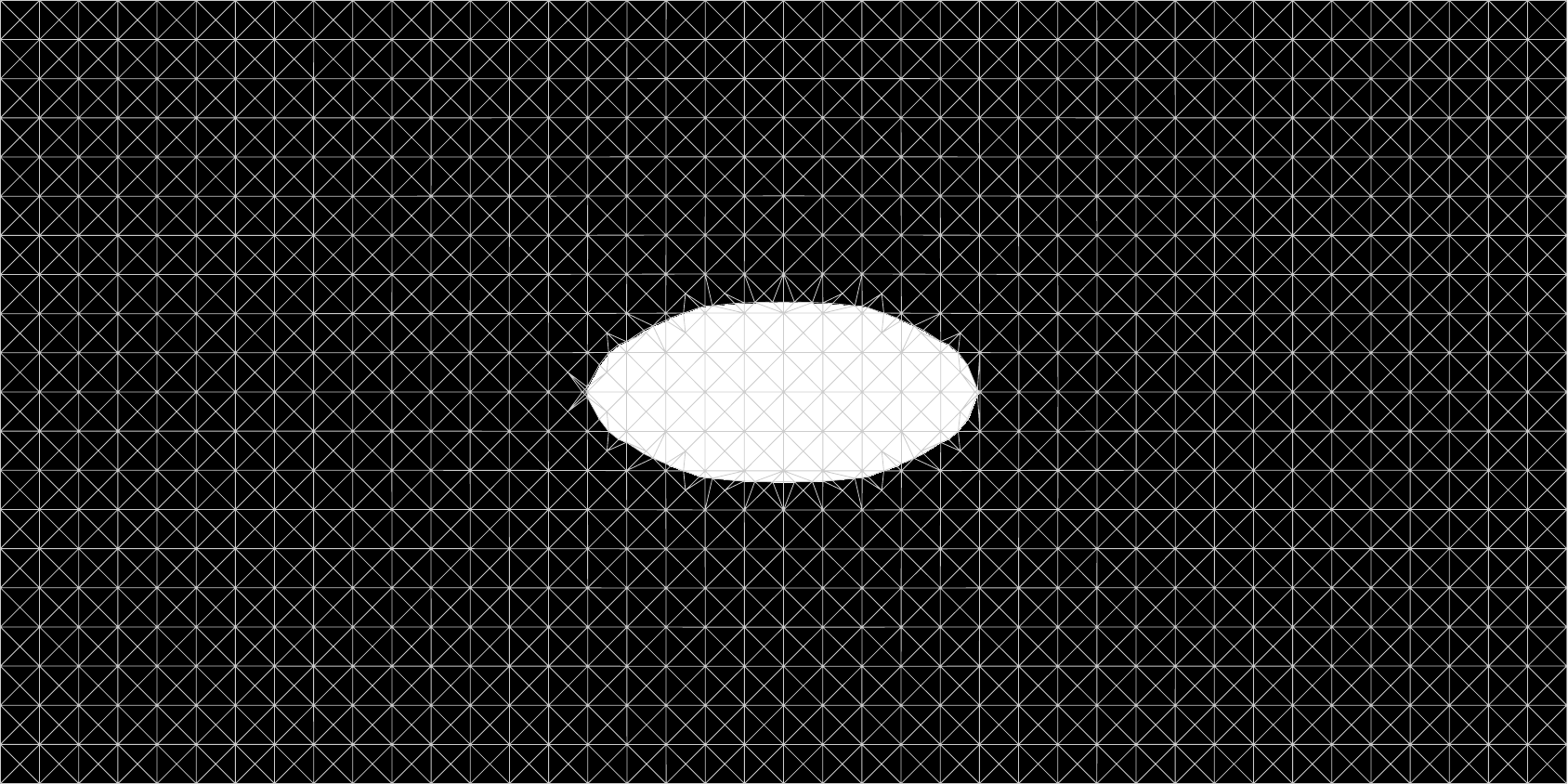}
        \caption{Initial domain for the sharp method (blow up).}
        \label{fig:navier_stokes:sharp:projected}
    \end{subfigure}
    \\
    \begin{subfigure}[t]{.45\linewidth}
        \centering
        \includegraphics[trim=180 100 180 100,clip,width = \linewidth]{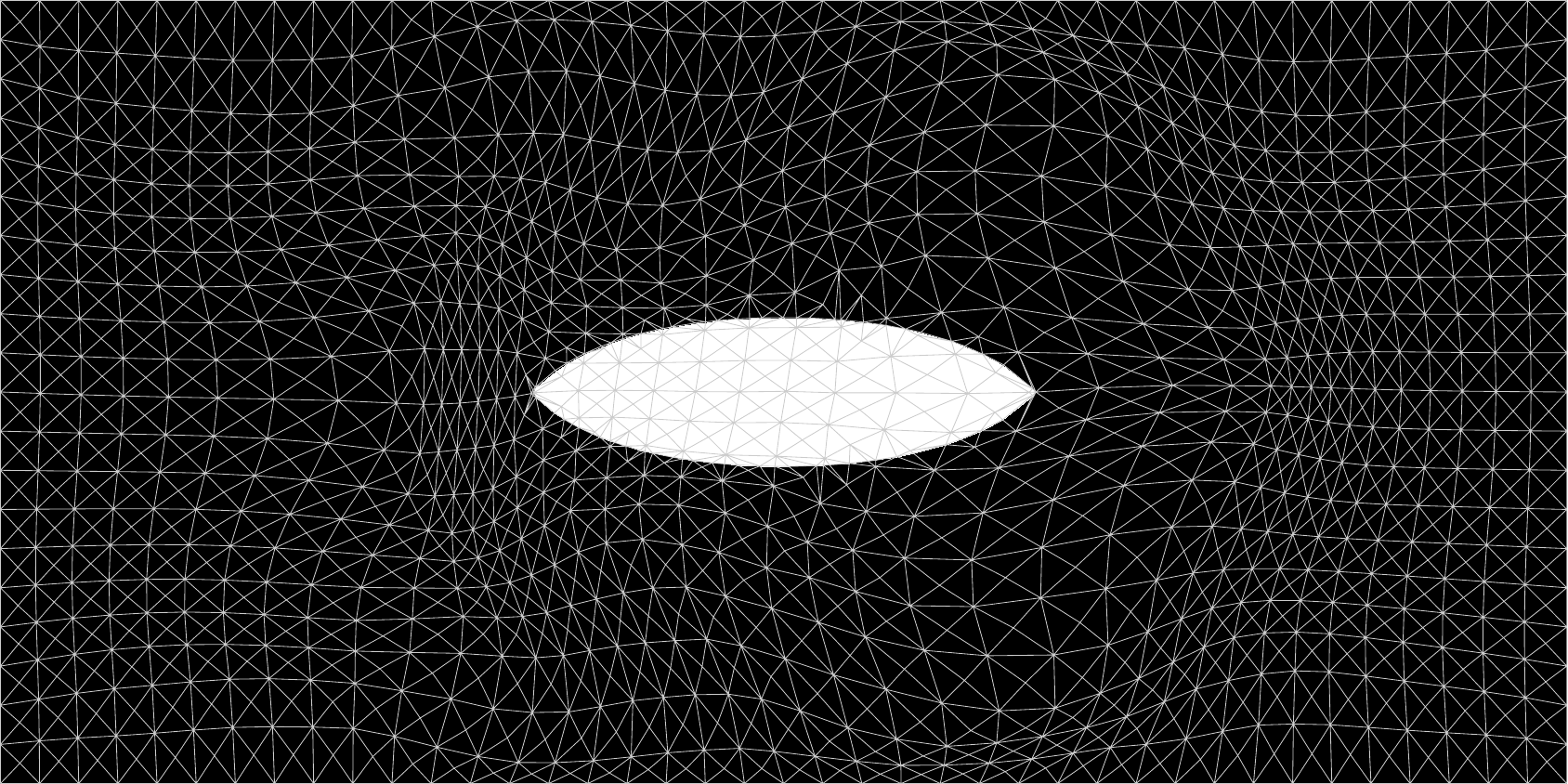}
        \caption{Final domain for the sharp method (blow up).}
        \label{fig:navier_stokes:sharp:final-blowup}
    \end{subfigure}
    \hfill    
    \begin{subfigure}[t]{.45\linewidth}
        \centering
        \includegraphics[width = \linewidth]{img/navierStokes/sharp/final/navierStokesSharpfinal_cell_mesh.pdf}
        \caption{Final domain for the sharp method.}
        \label{fig:navier_stokes:sharp:final}
    \end{subfigure}
    \caption{Numerical solution to problem \eqref{prob:NS:sharp}. Figure~\ref{fig:navier_stokes:sharp:final} shows the full domain $D$, the other figures show a blow-up of the relevant domain around the obstacle.
    Figure~\ref{fig:navier_stokes:diffuse:initial} and 
    Figure~\ref{fig:navier_stokes:diffuse:final} show the initial and optimised phase field function.
        In Figure~\ref{fig:navier_stokes:sharp:initial} we show the mesh after application of  Algorithm~\ref{alg:p2s:phase2sharp} and in Figure~\ref{fig:navier_stokes:sharp:projected} the projection of this, as described in Section~\ref{sec:Tech:geometricProjection}, which is the initial condition for the sharp setting.
        Figure~\ref{fig:navier_stokes:sharp:final-blowup} and Figure~\ref{fig:navier_stokes:sharp:final} both show the results of the $W^{1,\infty}$ shape optimisation after 20 steps.
The optimised shape has a larger aspect ratio than the initial guess of the sharp method, which is explained by the fact that the phase field optimisation problem contains perimeter regularisation and thus tends to find rounder solutions.}
    \label{fig:NS:results}
\end{figure}

We obtained the numerical results, that are depicted in Figure~\ref{fig:NS:results}. 
Since the relevant aspect of the domain is the fixed obstacle, we present magnified views of the numerical results.
The initial phase field is almost constant, with only minor variations, as shown in Figure~\ref{fig:navier_stokes:diffuse:initial}. 
Figure~\ref{fig:navier_stokes:diffuse:final} shows the optimized phase field in the shape of a rugby ball. After applying Algorithm~\ref{alg:p2s:phase2sharp}, this is depicted in Figure~\ref{fig:navier_stokes:sharp:initial}.
The projection of this domain onto one with the (up to numerical tolerance) correct volume and barycentre is shown in Figure \ref{fig:navier_stokes:sharp:projected}.
Finally, in Figure~\ref{fig:navier_stokes:sharp:final}, we depict the optimized shape, after a maximum of $20$ optimisation steps, in the form of an airfoil, as expected for this setup.
A comparison with the optimised phase field clearly shows the influence of the inherent surface area regularisation contained in the phase field approximation. 
The optimised shape is significantly pointier on the left and right than the optimised phase field.

In Figure~\ref{fig:navier_stokes:sharp:final} we also clearly observe the strength of the $W^{1,\infty}$ method in maintaining local mesh quality.

\section{Conclusion}

In shape optimization algorithms based on domain variations, good initial shapes that include the correct topology are crucial to maintaining mesh quality. 
We propose using phase field approximations of shape optimisation problems to generate these initial shapes. This approach ensures that the initial domains contain a more informed topology of the optimised shapes.

This requires a concept for generating an initial mesh for the shape optimisation method based on an implicitly given shape from a phase-field formulation. We propose updating the phase-field mesh by conformingly adding the zero-level line of the phase field, as well as post-processing the newly added vertices to maintain mesh quality in terms of regularity.  

The initial shape and mesh are then extracted from this based on the signum of the phase field at cell-barycentres. 
This post processing step may be exchanged for another method.

A numerical comparison with a simpler strategy for transferring information for the phase field solution indicates that inserting the zero-level line significantly improves both the optimised shape and the procedure.
We have also demonstrated that this approach enables shape optimisation methods to be used for structural optimisation, even when no a priori information on optimal topologies is available.

In future work, this approach could incorporate adaptivity.

\section*{Statements and Declarations}

\subsubsection*{Conflict of Interest.}
The authors declared that they have no conflict of interest.

\subsubsection*{Funding.}
M.H. and C.K. acknowledge funding of the project \textit{Fluiddynamische Formoptimierung mit Phasenfeldern und Lipschitz-Methoden} by the German Research Foundation under the project number \href{https://gepris.dfg.de/gepris/projekt/543959359?language=en}{543959359}.

\subsubsection*{Author Contribution.}
(CRediT taxonomy) 
\begin{description}
    \item[P.H.] Conceptualization, Investigation, Methodology, Project Administration, Software, Visualization, Writing -- original draft, Writing -- review \& editing
    \item[M.H.] Conceptualization, Funding Acquisition, Investigation, Methodology, Project Administration,  Writing -- original draft, Writing -- review \& editing
    \item[C.K.] Conceptualization, Funding Acquisition, Investigation, Methodology, Project Administration,  Writing -- original draft, Writing -- review \& editing
\end{description}

\subsubsection*{Acknowledgement.}
M.H. and C.K. acknowledge funding of the project \textit{Fluiddynamische Formoptimierung mit Phasenfeldern und Lipschitz-Methoden} by the German Research Foundation under the project number \href{https://gepris.dfg.de/gepris/projekt/543959359?language=en}{543959359}.

\subsubsection*{Usage of artificial intelligence.}
The software \texttt{DeepL Write} was used to improve language quality of the original draft.


\begingroup
\printbibliography

@book{henrot2018shape,
  title={Shape variation and optimization},
  author={Henrot, Antoine and Pierre, Michel},
  year={2018},
  publisher={European Mathematical Society-EMS-Publishing House GmbH}
}

@book{HPUU2009,
	title={Optimization with PDE Constraints},
	author={M. Hinze and R. Pinnau and M. Ulbrich and S. Ulbrich},
	series={Mathematical Modelling: Theory and Applications},
	year={2009},
	publisher={Springer Dordrecht},
doi={10.1007/978-1-4020-8839-1},
}

@book{BucBut05,
  title={Variational methods in shape optimization problems},
  author={Bucur, Dorin and Buttazzo, Giuseppe},
  year={2005},
  publisher={Birkhäuser Boston},
    doi={10.1007/b137163}
}

@book{DZ11,
  title={Shapes and geometries: metrics, analysis, differential calculus, and optimization},
  author={Delfour, Michel C and Zol{\'e}sio, J-P},
  year={2011},
  publisher={SIAM},
doi={10.1137/1.9780898719826},
}

@article{barrett1987fitted,
  title={Fitted and unfitted finite-element methods for elliptic equations with smooth interfaces},
  author={Barrett, John W and Elliott, Charles M},
  journal={IMA journal of numerical analysis},
  volume={7},
  number={3},
  pages={283--300},
  year={1987},
  publisher={Oxford University Press},
  ISSN = {1464-3642},
  url = {http://dx.doi.org/10.1093/imanum/7.3.283},
  DOI = {10.1093/imanum/7.3.283},
}

@article{wang2022sequentially,
  title={Sequentially coupled gradient-based topology and domain shape optimization},
  author={Wang, Zhijun and Suiker, Akke SJ and Hofmeyer, H{\`e}rm and Kalkman, Ivo and Blocken, Bert},
  journal={Optimization and Engineering},
  volume={23},
  number={1},
  pages={25--58},
  year={2022},
  publisher={Springer},
  ISSN = {1573-2924},
  url = {http://dx.doi.org/10.1007/s11081-020-09546-3},
  DOI = {10.1007/s11081-020-09546-3},
}

@article{auricchio2024analysis,
  title = {Analysis of a Combined Filtered/Phase-Field Approach to Topology Optimization in Elasticity},
  volume = {89},
  ISSN = {1432-0606},
  url = {http://dx.doi.org/10.1007/s00245-024-10104-x},
  DOI = {10.1007/s00245-024-10104-x},
  number = {2},
  journal = {Applied Mathematics \& Optimization},
  publisher = {Springer Science and Business Media LLC},
  author = {Auricchio,  Ferdinando and Marino,  Michele and Mazari,  Idriss and Stefanelli,  Ulisse},
  year = {2024},
  month = {02} 
}

@article{wang2007extended,
title = {An extended level set method for shape and topology optimization},
journal = {Journal of Computational Physics},
volume = {221},
number = {1},
pages = {395-421},
year = {2007},
issn = {0021-9991},
doi = {10.1016/j.jcp.2006.06.029},
url = {https://www.sciencedirect.com/science/article/pii/S0021999106002968},
author = {S.Y. Wang and K.M. Lim and B.C. Khoo and M.Y. Wang},
}

@article{starke2024shape,
author = {Starke, Gerhard},
title = {Shape Optimization by Constrained First-Order System Least Mean Approximation},
journal = {SIAM Journal on Scientific Computing},
volume = {46},
number = {5},
pages = {A3044-A3066},
year = {2024},
doi = {10.1137/23M1605570},
URL = {https://doi.org/10.1137/23M1605570},
eprint = {https://doi.org/10.1137/23M1605570}
}

@article{luft2021simultaneous,
  title={Simultaneous shape and mesh quality optimization using pre-shape calculus},
  author={Luft, Daniel and Schulz, Volker},
  journal={Control and Cybernetics},
  volume={50},
  number={4},
  pages={473--520},
  year={2021},
  publisher={Polska Akademia Nauk. Instytut Bada{\'n} Systemowych PAN}
}

@article{chen2026inertial,
  title={An inertial minimal-deformation-rate framework for shape optimization},
  author={Chen, Falai and Li, Buyang and Li, Jiajie and Tang, Rong},
  journal={arXiv preprint arXiv:2601.22605},
doi = {10.48550/arXiv.2601.22605},
  year={2026}
}

@article{blauth2026enforcing,
  title={Enforcing mesh quality constraints in shape optimization with a gradient projection method},
  author={Blauth, Sebastian and Leith{\"a}user, Christian},
  journal={Computer Methods in Applied Mechanics and Engineering},
  volume={448},
  pages={118451},
  year={2026},
  publisher={Elsevier},
  ISSN = {0045-7825},
  url = {http://dx.doi.org/10.1016/j.cma.2025.118451},
  DOI = {10.1016/j.cma.2025.118451},
}

@article{iglesias2018two,
author = {Iglesias, Jos\'{e} A. and Sturm, Kevin and Wechsung, Florian},
title = {Two-Dimensional Shape Optimization with Nearly Conformal Transformations},
journal = {SIAM Journal on Scientific Computing},
volume = {40},
number = {6},
pages = {A3807-A3830},
year = {2018},
doi = {10.1137/17M1152711},
URL = {https://doi.org/10.1137/17M1152711},
eprint = {https://doi.org/10.1137/17M1152711},
}

@software{DunePythonReference,
  author       = {Dedner, Andreas and
                  Kloefkorn, Robert and
                  Nolte, Martin},
  title        = {Python Bindings for the DUNE-FEM module},
  month        = mar,
  year         = 2020,
  publisher    = {Zenodo},
  version      = {v2.7.0},
  doi          = {10.5281/zenodo.3706994},
  url          = {https://doi.org/10.5281/zenodo.3706994}
}

@book{DUNE-book,
  title = {DUNE — The Distributed and Unified Numerics Environment},
  ISBN = {9783030597023},
  ISSN = {2197-7100},
  url = {http://dx.doi.org/10.1007/978-3-030-59702-3},
  DOI = {10.1007/978-3-030-59702-3},
  journal = {Lecture Notes in Computational Science and Engineering},
  publisher = {Springer International Publishing},
  author = {Sander,  Oliver},
  year = {2020}
}

@article{DecHerHin22,
	author = {Deckelnick, Klaus and Herbert, Philip J. and Hinze, Michael},
	title = {A novel $W^{1,\infty}$ approach to shape optimisation with Lipschitz domains},
	DOI= "10.1051/cocv/2021108",
	url= "https://doi.org/10.1051/cocv/2021108",
	journal = {ESAIM: COCV},
	year = 2022,
	volume = 28,
	pages = "2",
}

@article{DecHerHin24-App,
    author = {Klaus Deckelnick and Philip J. Herbert and Michael Hinze},
    title = {PDE constrained shape optimisation with first-order and Newton-type methods in the $W^{1,\infty}$ topology},
      ISSN = {1029-4937},
      url = {http://dx.doi.org/10.1080/10556788.2024.2424525},
      DOI = {10.1080/10556788.2024.2424525},
      journal = {Optimization Methods and Software},
      publisher = {Informa UK Limited},
      author = {Deckelnick,  Klaus and Herbert,  Philip J. and Hinze,  Michael},
      year = {2024},
      month = dec,
      pages = {1–27}
}

@article{DecHerHin25-NA,
  title = {Convergence of a steepest descent algorithm in shape optimisation using $W^{1,\infty}$ functions},
  volume = {59},
  ISSN = {2804-7214},
  url = {http://dx.doi.org/10.1051/m2an/2025033},
  DOI = {10.1051/m2an/2025033},
  number = {3},
  journal = {ESAIM: Mathematical Modelling and Numerical Analysis},
  publisher = {EDP Sciences},
  author = {Deckelnick,  Klaus and Herbert,  Philip Justin and Hinze,  Michael},
  year = {2025},
  month = {05},
  pages = {1505–1529}
}

@article{MulPinRun23,
author = {M\"{u}ller, Peter Marvin and Pinz\'{o}n, Jos\'{e} and Rung, Thomas and Siebenborn, Martin},
title = {A Scalable Algorithm for Shape Optimization with Geometric Constraints in Banach Spaces},
journal = {SIAM Journal on Scientific Computing},
volume = {45},
number = {2},
pages = {B231-B251},
year = {2023},
doi = {10.1137/22M1494609},
URL = {https://doi.org/10.1137/22M1494609},
eprint = {https://doi.org/10.1137/22M1494609}
}

@article{HerPinSie23,
  title = {Shape Optimization in $W^{1, \infty }$ with Geometric Constraints: a Study in Distributed-Memory Systems},
  volume = {207},
  ISSN = {1573-2878},
  url = {http://dx.doi.org/10.1007/s10957-025-02802-5},
  DOI = {10.1007/s10957-025-02802-5},
  number = {43},
  journal = {Journal of Optimization Theory and Applications},
  publisher = {Springer Science and Business Media LLC},
  author = {Herbert,  Philip J. and Escobar,  Jose A. Pinzon and Siebenborn,  Martin},
  year = {2025},
  month = {08} 
}

@article{MulKuhSie21,
  title={A novel p-harmonic descent approach applied to fluid dynamic shape optimization},
  author={M{\"u}ller, Peter Marvin and K{\"u}hl, Niklas and Siebenborn, Martin and Deckelnick, Klaus and Hinze, Michael and Rung, Thomas},
  journal={Structural and multidisciplinary optimization},
  volume={64},
  number={6},
  pages={3489--3503},
  year={2021},
  publisher={Springer},
ISSN = {1615-1488},
  url = {http://dx.doi.org/10.1007/s00158-021-03030-x},
  DOI = {10.1007/s00158-021-03030-x},
}

@article{LufSch21-A,
  title={Pre-shape calculus and its application to mesh quality optimization},
  author={Luft, Daniel and Schulz, Volker},
  journal={Control and Cybernetics},
  volume={50},
  number={3},
  pages={263--301},
  year={2021}
}

@article{HerLoa23,
  title={A manifold of planar triangular meshes with complete Riemannian metric},
  author={Herzog, Roland and Loayza-Romero, Estefan{\'\i}a},
  journal={Mathematics of Computation},
  volume={92},
  number={339},
  pages={1--50},
  year={2023},
doi          = {10.1090/MCOM/3775},
}

@inbook{AllDapJou21,
  title={Shape and topology optimization},
  author={Allaire, Gr{\'e}goire and Dapogny, Charles and Jouve, Fran{\c{c}}ois},
  booktitle={Handbook of numerical analysis},
  volume={22},
  pages={1--132},
  year={2021},
  publisher={Elsevier},
  ISBN = {9780444643056},
  ISSN = {1570-8659},
  url = {http://dx.doi.org/10.1016/bs.hna.2020.10.004},
  DOI = {10.1016/bs.hna.2020.10.004},
}

@article{PapFarSur-2021-MultipleSolutionsTopoOpt,
author = {Papadopoulos, Ioannis P. A. and Farrell, Patrick E. and Surowiec, Thomas M.},
title = {Computing Multiple Solutions of Topology Optimization Problems},
journal = {SIAM Journal on Scientific Computing},
volume = {43},
number = {3},
pages = {A1555-A1582},
year = {2021},
doi = {10.1137/20M1326209},
URL = {https://doi.org/10.1137/20M1326209},
arxiv = {2004.11797v2}
}

@article{EltHerLoa20,
  title={First and second order shape optimization based on restricted mesh deformations},
  author={Etling, Tommy and Herzog, Roland and Loayza, Estefan{\'\i}a and Wachsmuth, Gerd},
  journal={SIAM Journal on Scientific Computing},
  volume={42},
  number={2},
  pages={A1200--A1225},
  year={2020},
  publisher={SIAM},
  ISSN = {1095-7197},
  url = {http://dx.doi.org/10.1137/19M1241465},
  DOI = {10.1137/19m1241465},
}

@article {Modica,
	AUTHOR = {Modica, Luciano},
	TITLE = {The gradient theory of phase transitions and the minimal
	interface criterion},
	JOURNAL = {Arch. Rational Mech. Anal.},
	FJOURNAL = {Archive for Rational Mechanics and Analysis},
	VOLUME = {98},
	YEAR = {1987},
	NUMBER = {2},
	PAGES = {123--142},
	ISSN = {0003-9527},
	MRCLASS = {76T05 (80A15)},
	MRNUMBER = {866718},
	MRREVIEWER = {L. Hsiao},
	DOI = {10.1007/BF00251230},
	URL = {https://doi.org/10.1007/BF00251230},
}

@article {BlankR-VMPT,
    AUTHOR = {Blank, Luise and Rupprecht, Christoph},
     TITLE = {An extension of the projected gradient method to a {B}anach
              space setting with application in structural topology
              optimization},
   JOURNAL = {SIAM J. Control Optim.},
  FJOURNAL = {SIAM Journal on Control and Optimization},
    VOLUME = {55},
      YEAR = {2017},
    NUMBER = {3},
     PAGES = {1481--1499},
      ISSN = {0363-0129},
   MRCLASS = {49M05 (49M15 65K10 65K15 74P05 90C30 90C48)},
  MRNUMBER = {3648069},
       DOI = {10.1137/16M1092301},
       URL = {https://doi.org/10.1137/16M1092301},
}

@article{GarHHk-2015-NumApproxPhaseTopoOptFluids,
  title={Numerical {A}pproximation of {P}hase {F}ield {B}ased {S}hape and {T}opology {O}ptimization for {F}luids},
  author={Garcke, H. and Hecht, C. and Hinze, M. and Kahle, C.},
  journal={SIAM Journal on Scientific Computing},
  volume={37},
  number={4},
  pages={A1846--A1871},
  year={2015},
  publisher={Society for Industrial and Applied Mathematics},
  doi={10.1137/140969269},

}

@article{GarHKL-2018-TopOptPF-StateConstraints,
  title={A {P}hase {F}ield {A}pproach to {S}hape {O}ptimization in {N}avier--{S}tokes {F}low with {I}ntegral {S}tate {C}onstraints},
  author={Garcke, H. and Hinze, M. and Kahle, C. and Lam, K.F.},
  journal={Advances in Computational Mathematics},
  volume={44},
  number={5},
  pages={1345--1383},
  year={2018},
  publisher={Springer},
  doi={10.1007/s10444-018-9586-8},
}

@article{GarckeHKKL-2023-PhaseFieldShapeOptLaplace,
	author = {Garcke, Harald and H\"uttl, Paul and Kahle, Christian and Knopf, Patrik and Laux, Tim},
	title = {Phase-field methods for spectral shape and topology optimization},
	DOI= "10.1051/cocv/2022090",
	url= "https://doi.org/10.1051/cocv/2022090",
	journal = {ESAIM:  Control, Optimisation and Calculus of Variations},
	year = 2023,
	volume = 29,
	pages = "10",

}

@Article{GarckeHKK-2024-MultiPhaseElasticSpectral,
author={Garcke, Harald
and H{\"u}ttl, Paul
and Kahle, Christian
and Knopf, Patrik},
title={Sharp-Interface Limit of a Multi-phase Spectral Shape Optimization Problem for Elastic Structures},
journal={Applied Mathematics {\&} Optimization},
year={2024},
month={01},
day={09},
volume={89},
number={1},
pages={24},
issn={1432-0606},
doi={10.1007/s00245-023-10093-3},
url={https://doi.org/10.1007/s00245-023-10093-3}
}

@article{GarckeHuettlKnopf-2022-TopoOptEigValElasticity,
author = {Harald Garcke and Paul Hüttl and Patrik Knopf},
doi = {doi:10.1515/anona-2020-0183},
url = {https://doi.org/10.1515/anona-2020-0183},
title = {Shape and topology optimization involving the eigenvalues of an elastic structure: A multi-phase-field approach},
journal = {Advances in Nonlinear Analysis},
number = {1},
volume = {11},
year = {2022},
pages = {159--197},
}

@article{GarckeLamNuernbergSignori-2023-OverhangPenalizationManufacturingTopoOpt,
  title = {Overhang Penalization in Additive Manufacturing via Phase Field Structural Topology Optimization with Anisotropic Energies},
  volume = {87},
  ISSN = {1432-0606},
  url = {http://dx.doi.org/10.1007/s00245-022-09939-z},
  DOI = {10.1007/s00245-022-09939-z},
  number = {3},
  journal = {Applied Mathematics {\&} Optimization},
  publisher = {Springer Science and Business Media LLC},
  author = {Garcke,  Harald and Lam,  Kei Fong and N\"{u}rnberg,  Robert and Signori,  Andrea},
  year = {2023},
}

@Inbook{BlankEtAl-2012-PhaseFieldStructuralTopoOpt,
author="Blank, Luise
and Garcke, Harald
and Sarbu, Lavinia
and Srisupattarawanit, Tarin
and Styles, Vanessa
and Voigt, Axel",
editor="Leugering, G{\"u}nter
and Engell, Sebastian
and Griewank, Andreas
and Hinze, Michael
and Rannacher, Rolf
and Schulz, Volker
and Ulbrich, Michael
and Ulbrich, Stefan",
title="Phase-field Approaches to Structural Topology Optimization",
bookTitle="Constrained Optimization and Optimal Control for Partial Differential Equations",
year="2012",
publisher="Springer Basel",
address="Basel",
pages="245--256",
isbn="978-3-0348-0133-1",
doi="10.1007/978-3-0348-0133-1_13",
url="https://doi.org/10.1007/978-3-0348-0133-1_13"
}

@Article{DedeBordenHughes-2012-IsogeometricalTopologyOpt,
author={Ded{\`e}, Luca and Borden, Micheal J. and Hughes, Thomas J. R.},
title={Isogeometric Analysis for Topology Optimization with a Phase Field Model},
journal={Archives of Computational Methods in Engineering},
year={2012},
day={01},
volume={19},
number={3},
pages={427-465},
issn={1886-1784},
doi={10.1007/s11831-012-9075-z},
url={https://doi.org/10.1007/s11831-012-9075-z}
}

@article{Bartels-2015-RobustnessErrorEstimatesPhaseFieldTopoChange,
  author = {Bartels, S{\"o}ren},
  title = {Robustness of error estimates for phase field models at a
              class of topological changes},
  journal = {Comput. Methods Appl. Mech. Engrg.},
  fjournal = {Computer Methods in Applied Mechanics and Engineering},
  volume = {288},
  year = {2015},
  pages = {75--82},
  issn = {0045-7825},
  mrclass = {35B25 (35K91 74A50 80A22)},
  mrnumber = {3327018},
  doi = {10.1016/j.cma.2014.11.005},
  url = {http://dx.doi.org/10.1016/j.cma.2014.11.005}
}

@Article{Evgrafov-2005-LimitsOfPorousMediumTopologyOptimization,
author={Evgrafov, Anton},
title={The Limits of Porous Materials in the Topology Optimization of Stokes Flows},
journal={Applied Mathematics and Optimization},
year={2005},
day={01},
volume={52},
number={3},
pages={263-277},
issn={1432-0606},
doi={10.1007/s00245-005-0828-z},
url={https://doi.org/10.1007/s00245-005-0828-z}
}

@article{TakezawaNishiwakiKitamura-2010-ShapeTopoOptPhaseFieldSensitivity,
title = {Shape and topology optimization based on the phase field method and sensitivity analysis},
journal = {Journal of Computational Physics},
volume = {229},
number = {7},
pages = {2697-2718},
year = {2010},
issn = {0021-9991},
doi = {10.1016/j.jcp.2009.12.017},
author = {Akihiro Takezawa and Shinji Nishiwaki and Mitsuru Kitamura},
}

@Article{ZhouWang-2006-MulitmaterialStructuralTopoOptGeneralizedCHMultiphase,
author={Zhou, Shiwei
and Wang, Michael Yu},
title={Multimaterial structural topology optimization with a generalized Cahn--Hilliard model of multiphase transition},
journal={Structural and Multidisciplinary Optimization},
year={2006},
day={18},
volume={33},
number={2},
pages={89},
issn={1615-1488},
doi={10.1007/s00158-006-0035-9},
}

@Article{Wallin2012,
author={Wallin, Mathias
and Ristinmaa, Matti
and Askfelt, Henrik},
title={Optimal topologies derived from a phase-field method},
journal={Structural and Multidisciplinary Optimization},
year={2012},
day={01},
volume={45},
number={2},
pages={171-183},
doi={10.1007/s00158-011-0688-x},

}

@article {AlmiStefanelli,
	AUTHOR = {Almi, Stefano and Stefanelli, Ulisse},
	TITLE = {Topology optimization for incremental elastoplasticity: a
	phase-field approach},
	JOURNAL = {SIAM J. Control Optim.},
	FJOURNAL = {SIAM Journal on Control and Optimization},
	VOLUME = {59},
	YEAR = {2021},
	NUMBER = {1},
	PAGES = {339--364},
	ISSN = {0363-0129},
	MRCLASS = {74C05 (49J20 49K20 49Q10 74P10)},
	MRNUMBER = {4202025},
	DOI = {10.1137/20M1331275},
	URL = {https://doi.org/10.1137/20M1331275},
}

@article{Auricchio2,
author = {Auricchio, Ferdinando and Bonetti, Elena and Carraturo, Massimo and H\"{o}mberg, Dietmar and Reali, Alessandro and Rocca, Elisabetta},
title = {A phase-field-based graded-material topology optimization with stress constraint},
journal = {Mathematical Models and Methods in Applied Sciences},
volume = {30},
number = {08},
pages = {1461-1483},
year = {2020},
doi = {10.1142/S0218202520500281},
URL = {https://doi.org/10.1142/S0218202520500281},
eprint = {https://doi.org/10.1142/S0218202520500281}
}

@InProceedings{Bourdin-Chambolle,
	author="Bourdin, B.	and Chambolle, A.",
	editor="Bends{\o}e, M.P. and Olhoff, N. and Sigmund, O.",
	title="{The Phase-Field Method in Optimal Design}",
	booktitle="IUTAM Symposium on Topological Design Optimization of Structures, Machines and Materials",
	year="2006",
	publisher="Springer Netherlands",
	address="Dordrecht",
    pages="207--215",
	isbn="978-1-4020-4752-7",
    doi = "10.1007/1-4020-4752-5_21",
url = "https://doi.org/10.1007/1-4020-4752-5_21"
}

@article {Burger-Stainko,
	AUTHOR = {Burger, M. and Stainko, R.},
	TITLE = {Phase-field relaxation of topology optimization with local
	stress constraints},
	JOURNAL = {SIAM J. Control Optim.},
	FJOURNAL = {SIAM Journal on Control and Optimization},
	VOLUME = {45},
	YEAR = {2006},
	NUMBER = {4},
	PAGES = {1447--1466},
	ISSN = {0363-0129},
	MRCLASS = {74P15 (74P05 74P10 74S05 90C51)},
	MRNUMBER = {2257229},
	MRREVIEWER = {Michal Ko\v{c}vara},
	DOI = {10.1137/05062723X},
	URL = {https://doi.org/10.1137/05062723X},
}

@article {Carraturo,
    AUTHOR = {Carraturo, Massimo and Rocca, Elisabetta and Bonetti, Elena
              and H\"{o}mberg, Dietmar and Reali, Alessandro and Auricchio,
              Ferdinando},
     TITLE = {Graded-material design based on phase-field and topology
              optimization},
   JOURNAL = {Comput. Mech.},
  FJOURNAL = {Computational Mechanics},
    VOLUME = {64},
      YEAR = {2019},
    NUMBER = {6},
     PAGES = {1589--1600},
      ISSN = {0178-7675},
   MRCLASS = {74P15 (74P05)},
  MRNUMBER = {4031822},
MRREVIEWER = {Bernard Rousselet},
       DOI = {10.1007/s00466-019-01736-w},
       URL = {https://doi.org/10.1007/s00466-019-01736-w},
}

@article {Penzler,
	AUTHOR = {Penzler, P. and Rumpf, M. and Wirth, B.},
	TITLE = {A phase-field model for compliance shape optimization in
	nonlinear elasticity},
	JOURNAL = {ESAIM Control Optim. Calc. Var.},
	FJOURNAL = {ESAIM. Control, Optimisation and Calculus of Variations},
	VOLUME = {18},
	YEAR = {2012},
	NUMBER = {1},
	PAGES = {229--258},
	ISSN = {1292-8119},
	MRCLASS = {49Q10 (74B20 74P05)},
	MRNUMBER = {2887934},
	MRREVIEWER = {Shiah-Sen Wang},
	DOI = {10.1051/cocv/2010045},
	URL = {https://doi.org/10.1051/cocv/2010045},
}
\endgroup

\end{document}